\documentclass[a4paper]{amsart}
\usepackage{amscd, amssymb, color, booktabs}
\usepackage{myown}
\title
{Generalized Preservation Principle in Finite Theta Correspondence}
\author{Shu-Yen Pan}
\address{Department of Mathematics,
National Tsing Hua University, Hsinchu 300, Taiwan}
\email{sypan@math.nthu.edu.tw}

\thanks{This work is partially supported by Taiwan MOST-grant 106-2115-M-007-015.}


\date{\today}

\begin{document}

\begin{abstract}
It is known that irreducible cuspidal characters satisfy the preservation principle
in the Howe correspondences of finite reductive dual pairs.
In this article, we generalize the preservation principle to any irreducible characters
of classical groups.
\end{abstract}

\maketitle
\tableofcontents

\section{Introduction}

\subsection{}\label{0102}
Let $(\bfG,\bfG')$ be a reductive dual pair over a finite field $\bfF_q$ of odd characteristic.
By restricting the Weil character (\cf.~\cite{gerardin}) to the dual pair $(\bfG,\bfG')$ with respect to
a nontrivial additive character $\psi$ of $\bfF_q$,
we have the \emph{Weil character} $\omega^\psi_{\bfG,\bfG'}$ of $G\times G'$,
which has non-negative integral decomposition
\[
\omega^\psi_{\bfG,\bfG'}=\sum_{\rho\in\cale(G),\ \rho'\in\cale(G')}m_{\rho,\rho'}\rho\otimes\rho'
\]
where $\cale(G)$ denotes the set of irreducible characters of group of rational points $G$ of $\bfG$.
We define the \emph{Howe correspondence} 
\[
\Theta^\psi_{\bfG,\bfG'}=\{\,(\rho,\rho')\in\cale(G)\times\cale(G')\mid m_{\rho,\rho'}\neq 0\,\}
\]
which is a relation between $\cale(G)$ and $\cale(G')$.
We say that an irreducible character $\rho\in\cale(G)$ \emph{occurs} in the correspondence 
for the dual pair $(\bfG,\bfG')$ if there exists $\rho'\in\cale(G')$ such that 
$(\rho,\rho')\in\Theta^\psi_{\bfG,\bfG'}$.
 
Now we consider the following types of \emph{related Witt series}:
\begin{enumerate}
\item[(I)] $\bfv$ is a hermitian space (over a quadratic extension of $\bfF_q$),
and $\{\bfv'^+\}$ (resp.~$\{\bfv'^-\}$) is the Witt series of even-dimensional
(resp.~odd-dimensional) hermitian spaces;

\item[(II)] $\bfv$ is a quadratic space,
and $\{\bfv'\}$ is the Witt series of symplectic spaces;

\item[(III)] $\bfv$ is a symplectic space,
and $\{\bfv'^+\}$ (resp.~$\{\bfv'^-\}$) is the Witt series of even-dimensional quadratic spaces
with trivial (resp.~two-dimensional) anisotropic kernel;

\item[(IV)] $\bfv$ is a symplectic space,
and $\{\bfv'\}$ is the Witt series of odd-dimensional quadratic spaces.
\end{enumerate}
For each case the pair $(G(\bfv),G(\bfv'))$ (or $(G(\bfv),G(\bfv'^\epsilon))$ where $\epsilon=\pm$)
forms a reductive dual pair over $\bfF_q$ where $G(\bfv)$ denotes the group of isometries of $\bfv$
with the given form.

Let $\rho\in\cale(G(\bfv))$.
For (I) and (III),
let $n_0'^\epsilon(\rho)$ denote the minimal dimension of $\bfv'^\epsilon$ such that
$\rho$ occurs in the Howe correspondence for the dual pair $(G(\bfv),G(\bfv'^\epsilon))$.
For (II) and (IV),
let $n_0'(\rho)$ denote the minimal dimension of $\bfv'$ such that
$\rho$ occurs in the Howe correspondence for the dual pair $(G(\bfv),G(\bfv'))$.
The following \emph{preservation principle} for cuspidal characters is from \cite{pn-depth0} theorem 12.3:

\begin{prop}\label{0103}
Let $\rho$ be an irreducible cuspidal character of $G(\bfv)$.
\begin{enumerate}
\item[(i)]
If $G(\bfv)=\rmU(\bfv)$, then
$n'^+_0(\rho)+n'^-_0(\rho)=2\,\dim(\bfv)+1$.

\item[(ii)] If $G(\bfv)=\rmO(\bfv)$, then
$n_0'(\rho)+n_0'(\rho\cdot\sgn)=2\,\dim(\bfv)$.

\item[(iii)] If $G(\bfv)=\Sp(\bfv)$ and $\bfv'^\pm$ even-dimensional, then
$n'^+_0(\rho)+n'^-_0(\rho)=2\,\dim(\bfv)+2$.

\item[(iv)] If $G(\bfv)=\Sp(\bfv)$ and $\bfv'$ odd-dimensional, then
$n'_0(\rho)+n'_0(\rho^c)=2\,\dim(\bfv)+2$.
\end{enumerate}
\end{prop}
Note that here ``$\sgn$'' denotes the sign character of\/ $\rmO(\bfv)$,
and ``$\rho^c$'' denotes the character obtained from $\rho$ by conjugating
certain element in the symplectic similitude group (\cf.~\cite{waldspurger} \S 4.11).

\subsection{}\label{0104}
The purpose of this article is to generalize the above preservation principle
to all irreducible characters.
To achieves this we need to recall Lusztig's parametrzation of irreducible
characters of finite classical groups.

Suppose that $\bfG=\rmU_n$, $\Sp_{2n}$, $\rmO^\epsilon_{2n}$, or $\SO_{2n+1}$.
Let $\cale(G)_s$ denote the \emph{Lusztig series} associated to the conjugacy class $(s)$
of a semisimple element $s$ in the connected component $(G^*)^0$ of the dual group $G^*$ of $G$.
In particular, $\cale(G)_1$ is the set of irreducible unipotent characters of $G$.
It is known that these Lusztig series partition the set of irreducible characters:
\[
\cale(G)=\bigcup_{(s)\subset (G^*)^0}\cale(G)_s.
\]
In \cite{lg} Lusztig shows that there exists a bijection
$\grL_s\colon\cale(G)_s\rightarrow\cale(C_{G^*}(s))_1$
where $C_{G^*}(s)$ denotes the centralizer of $s$.
Recall that we can define groups $G^{(j)}$ (where $j=0,1$ of $G$ is a unitary group,
and $j=0,1,2$ if $G$ is an orthogonal group or a symplectic group), and then
we have a bijection
\[
\Xi_s\colon\cale(G)_s\rightarrow\begin{cases}
\cale(G^{(0)}\times G^{(1)})_1, & \text{if $G$ is a unitary group};\\
\cale(G^{(0)}\times G^{(1)}\times G^{(2)})_1, & \text{otherwise}
\end{cases}
\]
modified from $\grL_s$ above.

A \emph{$\beta$-set} $A=\{a_1,a_2,\ldots,a_m\}$ is a finite subset of non-negative integers with
elements written in decreasing order, i.e., $a_1>a_2>\cdots>a_m$.
A \emph{symbol} (by Lusztig) $\Lambda=\binom{A}{B}$ is an ordered pair of two $\beta$-sets.
Lusztig gives a parametrization of the set $\cale(G)_1$ of irreducible unipotent characters of
a classical group $G$ in terms of certain set of symbols denoted by $\cals_\bfG$
(\cf.~Subsection~\ref{0234}, Subsection~\ref{0233}).
The unipotent character associated to $\Lambda$ is denoted by $\rho_\Lambda$.

For a $\beta$-set $A$ and a symbol $\Lambda=\binom{A}{B}$, we define
\begin{align}\label{0105}
\begin{split}
\delta(A) &=\begin{cases}
0, & \text{if $A=\emptyset$};\\
\max(A)-|A|+1, & \text{if $A\neq\emptyset$}
\end{cases} \\
\delta(\Lambda) &=\delta(A)+\delta(B)
\end{split}
\end{align}
where $\max(A)$ denotes the maximal element of $A$,
and $|A|$ denotes the number of elements in $A$.
For $\Lambda\in\cals_\bfG$, we define $\delta(\rho_\Lambda)=\rho(\Lambda)$.

Now we consider types (I), (II), (III), (IV) of related Witt series given in Subsection~\ref{0102}.
Let $\rho\in\cale(G(\bfv))$.
If $\bfv$ is not an odd-dimensional quadratic space,
we suppose that $\rho\in\cale(G(\bfv))_s$ for some $s$.
If $\bfv$ is an odd-dimensional quadratic space,
we let $G^0(\bfv)$ denote the connected component of $G(\bfv)$ and suppose that 
$\rho|_{G^0(\bfv)}\in\cale(G^0(\bfv))_s$ for some $s$.
Now we write
\[
\begin{cases}
\Xi_s(\rho)= \rho^{(0)}\otimes\rho^{(1)}, & \text{if $G(\bfv)$ is a unitary group};\\
\Xi_s(\rho|_{G^0(\bfv)})= \rho^{(0)}\otimes\rho^{(1)}\otimes\rho^{(2)}, & \text{if $G(\bfv)$ is an odd-orthogonal group}; \\
\Xi_s(\rho)= \rho^{(0)}\otimes\rho^{(1)}\otimes\rho^{(2)}, & \text{otherwise}
\end{cases} 
\]
(\cf.~(\ref{0235}), (\ref{0211})).
Then we define
\begin{equation}\label{0107}
\delta(\rho)=\begin{cases}
\delta(\rho^{(1)}), & \text{for types (I), (IV)};\\
\delta(\rho^{(2)}), & \text{for types (II), (III)}.
\end{cases}
\end{equation}
It is known that $\delta(\rho)=\delta(\rho^c)$ if $G$ is a symplectic group;
and $\delta(\rho)=\delta(\rho\cdot\sgn)$ if $G$ is an orthogonal group.

Now we state the generalization of the preservation principle
(\cf.~Theorems~\ref{0919}, \ref{0920}, \ref{0921}, and \ref{0702}) for any irreducible characters:

\begin{thm}\label{0106}
Let $\rho$ be an irreducible character of $G(\bfv)$.
\begin{enumerate}
\item[(i)]
If $G(\bfv)=\rmU(\bfv)$, then
$n'^+_0(\rho)+n'^-_0(\rho)=2\,\dim(\bfv)-2\delta(\rho)+1$.

\item[(ii)] If $G(\bfv)=\rmO(\bfv)$, then
$n_0'(\rho)+n_0'(\rho\cdot\sgn)=2\,\dim(\bfv)-2\delta(\rho)$.

\item[(iii)] If $G(\bfv)=\Sp(\bfv)$ and $\{\bfv'^\pm\}$ even-dimensional, then
$n'^+_0(\rho)+n'^-_0(\rho)=2\,\dim(\bfv)-2\delta(\rho)+2$.

\item[(iv)] If $G(\bfv)=\Sp(\bfv)$ and $\{\bfv'\}$ odd-dimensional, then
$n'_0(\rho)+n'_0(\rho^c)=2\,\dim(\bfv)-2\delta(\rho)+2$.
\end{enumerate}
\end{thm}
It is known that $\delta(\rho)=0$ if $\rho$ is cuspidal.
Hence, the above theorem is the generalization of the preservation principles
for cuspidal characters in Proposition~\ref{0103}.
Now we list two applications of our main theorem:
\begin{enumerate}
\item
For an irreducible character $\rho$ of a classical group $G$.
A non-negative integer called \emph{tensor rank} is attached to $\rho$ by Gurevich-Howe in
\cite{gurevich-howe-rank}.
The tensor rank is called \emph{$\Theta$-rank} and denoted by $\Theta\text{\rm -rk}(\rho)$ in \cite{pan-theta-rank}.
Then Theorem~\ref{0106} provides an upper bound for $\Theta\text{\rm -rk}(\rho)$.
More precisely, for $\rho\in\cale(G)$, we have
\[
\Theta\text{\rm -rk}(\rho)\leq\begin{cases}
2n, & \text{if $\bfG=\Sp_{2n},\rmO^\epsilon_{2n},\rmO_{2n+1}$};\\
n, & \text{if $\bfG=\rmU_n$}.
\end{cases}
\]
More details can be found in \cite{pan-theta-rank}.

\item In \cite{howe-finite} corollary 5.2.1, it is asserted that for each $0\leq k\leq n$, there exists
a unique $\rho_k\in\cale(\Sp_{2n}(q))$ such that $\rho_k$ occurs in the Howe correspondences for both
pairs $(\Sp_{2n},\rmO^+_{2k})$ and  $(\Sp_{2n},\rmO^-_{2(n-k)+2})$.
From Theorem~\ref{0106} we can see that $\rho_k$ is in fact unipotent and parametrized by the symbol
$\binom{n-k+1,0}{k}$ (\cf.~Remark~\ref{0301}).
\end{enumerate}

\subsection{}

Let $\bfG_n$ denote $\Sp_{2n}$, $\rmO_{2n+1}$, $\rmO^+_{2n}$, $\rmO^-_{2n+2}$, $\rmU_{2n}$ or $\rmU_{2n+1}$.
For $(\rho,\rho')\in\Theta^\psi_{\bfG_n,\bfG'_{n'}}$, we denote
\begin{align*}
\Theta_{\bfG'_{n'}}(\rho) &=\{\,\rho'\in\cale(\bfG'_{n'})\mid(\rho,\rho')\in\Theta^\psi_{\bfG_n,\bfG'_{n'}}\,\}, \\
\Theta_{\bfG_n}(\rho') &=\{\,\rho\in\cale(\bfG_n)\mid(\rho,\rho')\in\Theta^\psi_{\bfG_n,\bfG'_{n'}}\,\}.
\end{align*}
It is known that $\Theta_{\bfG'_{n'}}(\rho)\neq \emptyset$ implies that
$\Theta_{\bfG'_{n''}}(\rho)\neq\emptyset$ for any $n''\geq n'$.
Therefore, we say that $(\rho,\rho')\in\Theta^\psi_{\bfG_n,\bfG'_{n'}}$ is a \emph{first occurrence}
if $\Theta_{\bfG'_{n'-1}}(\rho)=\emptyset$ and $\Theta_{\bfG_{n-1}}(\rho')=\emptyset$.
It is well known that if $(\rho,\rho')$ occurs in the correspondence and both $\rho,\rho'$ are cuspidal,
then $(\rho,\rho')$ is a first occurrence.
Of course, the converse implication is in general not true.

Let $\rho\in\cale(G(\bfv))$.
It is shown in \cite{pan-Witt-series} that $\delta(\rho)=0$ if and only if there exists
$\rho'\in\cale(G(\bfv'))$ for some $\bfv'$ in its Witt series such that $(\rho,\rho')$ is a first occurrence.

\begin{cor}
Let $\rho$ be an irreducible character of a classical $G(\bfv)$.
Then $\rho$ satisfies the preservation principle in Proposition~\ref{0103}
if and only if there exists $\rho'\in\cale(G(\bfv'))$ for some $\bfv'$ in its Witt series
such that $(\rho,\rho')$ is a first occurrence in the Howe correspondence.
\end{cor}

\subsection{}\label{0101}
It is known that the unipotent characters are not preserved by the correspondence
for synmplectic/odd-orthogonal dual pairs.
To understand the correspondence for this kind of dual pairs, we need to introduce some new notations.
An irreducible character $\rho$ of $G$ is called \emph{pseudo-unipotent} if
$\langle\rho,R_{\bfT_w,\theta_w}\rangle\neq 0$ for some $\bfT_w$ where $\theta_w$ is an order $2$ character
of $\bfT_w$ (\cf.~\cite{pan-odd} \S 3.3).

When $G$ is a unitary group or an orthogonal group,
an irreducible character $\rho$ of $G$ is pseudo-unipotent if and only if
$\rho\chi_\bfG$ is unipotent
where $\chi_\bfG$ is a fixed character of $G$ of order $2$ (\cf.~Lemma~\ref{0405}).
However, the situation for symplectic groups is more interesting.

\begin{prop}\label{0404}
A symplectic group $\Sp_{2n}(q)$ has a pseudo-unipotent cuspidal character if and only if
$n=m^2$ for some integer $m$.
Moreover $\Sp_{2m^2}(q)$, $m>0$, has exactly two irreducible pseudo-unipotent cuspidal characters.
\end{prop}

The description of the first occurrences of unipotent or pseudo-unipotent cuspidal charatcers
for the dual pairs $(\Sp_{2n},\rmO^\epsilon_{2n'})$ or $(\Sp_{2n},\rmO_{2n'+1})$ are given
in Proposition~\ref{0409}, Proposition~\ref{0704}, and Proposition~\ref{0407}.

\subsection{}
The contents of this article are as follows.
In Section 2 we first give the definitions and notations which are used in this paper.
Then we recall some basic properties of Deligne-Lusztig virtual character and the characterization
of unipotent characters by Lusztig.
In Section 3 we discuss the ``\emph{generalized preservation principle}'' for the Howe correspondence
of irreducible unipotent characters.
In Section 4 we deduce the generalized preservation principle for general irreducible characters from
the preservation principle of irreducible unipotent characters discussed in Section 3.
In the final section, we define ``pseudo-unipotent characters'' for symplectic groups
and give a complete description of the Howe correspondence on unipotent and pseudo-unipotent
cuspidal characters.

When this article is nearly complete,
the author has a chance to read the preprint \cite{liu-wang} by Liu and Wang in which
some results are overlapped with the results in Section 5 of this article.
In particular, the ``pseudo-unipotent characters'' here are called ``$\theta$-representations''
in \cite{liu-wang}.
However, the approaches of two papers are somewhat different.

\section{Preliminaries}

\subsection{Lusztig's symbols}
Recall that a symbol $\Lambda=\binom{A}{B}$ is an ordered pair of two $\beta$-sets
(\cf.~Subsection~\ref{0104}).
The set $A$ (resp.~$B$) is called the \emph{first row} (resp.~the \emph{second row}) of $\Lambda$ and is also
denoted by $\Lambda^*$ (resp.~$\Lambda_*$).
One the set of symbols, we define an equivalence relation generated by
\[
\binom{a_1,a_2,\ldots,a_{m_1}}{b_1,b_2,\ldots,b_{m_2}}
\sim \binom{a_1+1,a_2+1,\ldots,a_{m_1}+1,0}{b_1+1,b_2+1,\ldots,b_{m_2}+1,0}.
\]
The set of equivalence classes of symbols is denoted by $\cals$.

For a symbol $\Lambda=\binom{a_1,\ldots,a_{m_1}}{b_1,\ldots,b_{m_2}}$,
we define its \emph{rank} and \emph{defect} by
\begin{align}\label{0213}
\begin{split}
{\rm rk}(\Lambda) &= \sum_{i=1}^{m_1}a_i+\sum_{j=1}^{m_2}b_j
-\biggl\lfloor\biggr(\frac{m_1+m_2-1}{2}\biggr)^2\biggr\rfloor, \\
{\rm def}(\Lambda) &=m_1-m_2.
\end{split}
\end{align}
It is easy to check that
\[
{\rm rk}(\Lambda)\geq\biggl\lfloor\biggl(\frac{{\rm def}(\Lambda)}{2}\biggr)^2\biggr\rfloor
\]
for any $\Lambda$.
A symbol  $\Lambda$ is called \emph{cuspidal} if the above inequality is in fact an equality.
It is known that $\Lambda$ is cuspidal if and only if the unipotent character $\rho_\Lambda$ is cuspidal.

Recall that a \emph{bi-partition} $\sqbinom{\mu}{\nu}$ of $n$ is an ordered pair of two
partitions $\mu,\nu$ such that $\bigl|\sqbinom{\mu}{\nu}\bigr|:=|\mu|+|\nu|=n$.
We define a mapping from symbols to bi-partitions by
\begin{equation}
\Upsilon\colon\binom{a_1,a_2,\cdots,a_{m_2}}{b_1,b_2,\cdots,b_{m_2}}\mapsto
\sqbinom{a_1-(m_1-1),a_2-(m_2-1),\ldots,a_{m_1}}{b_1-(m_2-1),b_2-(m_2-2),\ldots,b_{m_2}}.
\end{equation}
It is easy to see that $\Lambda_1\sim\Lambda_2$ if and only if $\Upsilon(\Lambda_1)=\Upsilon(\Lambda_2)$.
Moreover, it is not difficult to check that
\begin{equation}\label{0215}
{\rm rk}(\Lambda)
=|\Upsilon(\Lambda)|+\biggl\lfloor\biggl(\frac{{\rm def}(\Lambda)}{2}\biggr)^2\biggr\rfloor
=|\Upsilon(\Lambda)|+\begin{cases}
\bigl(\frac{{\rm def}(\Lambda)}{2}\bigr)^2, & \text{if ${\rm def}(\Lambda)$ even};\\
\bigl(\frac{{\rm def}(\Lambda)-1}{2}\bigr)\bigl(\frac{{\rm def}(\Lambda)+1}{2}\bigr),
& \text{if ${\rm def}(\Lambda)$ odd}.
\end{cases}
\end{equation}
Therefore $\Lambda$ is cuspidal if and only if $\Upsilon(\Lambda)$ is the
empty partition $\sqbinom{-}{-}$.

Recall that for a $\beta$-set $A$, a number $\delta(A)$ is given in Subsection~\ref{0104}.
It is clear that $\delta(A)\geq0$,
and $\delta(A)=0$ if and only if $A=\emptyset$ or $A=\{m-1,m-2,\ldots,0\}$ for some $m\geq 1$.
For a symbol $\Lambda=\binom{A}{B}$,
it is clear that $\delta(\Lambda)=\delta(\Lambda^\rmt)$ where $\binom{A}{B}^\rmt=\binom{B}{A}$,
and $\delta(\Lambda)=\delta(\Lambda')$ if $\Lambda\sim\Lambda'$.

\begin{lem}
A symbol $\Lambda$ is cuspidal if and only if $\delta(\Lambda)=0$.
\end{lem}
\begin{proof}
Suppose that $\delta(\Lambda)=0$ and write $\Lambda=\binom{A}{B}$.
Now $\delta(\Lambda)=0$ implies that $\delta(A)=\delta(B)=0$.
Hence $A$ is either empty or $A=\{m_1,m_1-1,\ldots,1,0\}$ for some $m_1\geq 0$.
Similarly, $B$ is either empty or $B=\{m_2,m_2-1,\ldots,1,0\}$ for some $m_2\geq 0$.
If $A$ or $B$ is empty, it is clear that $\Lambda$ is cuspidal.
If $A$ and $B$ are both non-empty,
without loss of generality, may assume that $m_1\geq m_2$,
then
\[
\Lambda=\binom{m_1,m_1-1,\ldots,1,0}{m_2,m_2-1,\ldots,1,0}
\sim\binom{m_1-m_2-1,m_1-m_2-2,\ldots,1,0}{-}.
\]
So $\Lambda$ is again cuspidal.
On the other hand,
it is also clear by definition that $\delta(\Lambda)=0$ if $\Lambda$ is cuspidal.
\end{proof}

\begin{lem}\label{0216}
Suppose that ${\rm rk}(\Lambda)=n$.
Then we have $\delta(\Lambda)\leq n$.
Moreover, $\delta(\Lambda)=n$ if and only if
\[
\Lambda\sim\begin{cases}
\binom{k}{n-k+1,0}, & \text{if\/ ${\rm def}(\Lambda)=-1$}; \\
\binom{n-k}{k}, & \text{if\/ ${\rm def}(\Lambda)=0$}; \\
\binom{n-k+1,0}{k}, & \text{if\/ ${\rm def}(\Lambda)=1$}
\end{cases}
\]
for some $0\leq k\leq n$.
\end{lem}
\begin{proof}
Suppose that ${\rm rk}(\Lambda)=n$ and write $\Lambda=\binom{a_1,a_2,\ldots,a_{m_1}}{b_1,b_2,\ldots,b_{m_2}}$.
Because $\delta(\Lambda)$ is invariant on the similar class of $\Lambda$,
we may assume that $m_1,m_2$ are both nonzero.
Now by definition we have
\[
\delta(\Lambda)=(a_1-m_1+1)+(b_1-m_2+1)\leq|\Upsilon(\Lambda)|
={\rm rk}(\Lambda)-\biggl\lfloor\biggl(\frac{{\rm def}(\Lambda)}{2}\biggr)^2\biggr\rfloor
\leq n.
\]
It is easy to check that if $\Lambda$ is similar to one of the symbols listed in the lemma,
then $\delta(\Lambda)=n$.
On the other hand, suppose that $\delta(\Lambda)=n$.
Then we must have $\bigl\lfloor\bigl(\frac{{\rm def}(\Lambda)}{2}\bigr)^2\bigr\rfloor=0$ and
$\Upsilon(\Lambda)=\sqbinom{a_1-m_1+1}{b_1-m_2+1}=\sqbinom{n-k}{k}$ for some $0\leq k\leq n$.
This implies that ${\rm def}(\Lambda)=0,\pm1$ and $\Lambda$ is similar to one the symbols listed in the lemma.
\end{proof}

\begin{lem}\label{0217}
Suppose that ${\rm rk}(\Lambda)=n$ and ${\rm def}(\Lambda)=\pm2$.
Then we have $\delta(\Lambda)\leq n-1$.
Moreover, $\delta(\Lambda)=n-1$ if and only if
\[
\Lambda\sim\textstyle\binom{k}{n-k+1,1,0},\binom{n-k+1,1,0}{k}
\]
for some $0\leq k\leq n-1$.
\end{lem}
\begin{proof}
Write $\Lambda=\binom{a_1,a_2,\ldots,a_{m_1}}{b_1,b_2,\ldots,b_{m_2}}$.
From the proof of the previous lemma, we see that
\[
\delta(\Lambda)
\leq |\Upsilon(\Lambda)|
= {\rm rk}(\Lambda)-\biggl\lfloor\biggl(\frac{{\rm def}(\Lambda)}{2}\biggr)^2\biggr\rfloor
=n-1.
\]
Moreover, if $\delta(\Lambda)=n-1$, we need
$\Upsilon(\Lambda)=\sqbinom{a_1-m_1+1}{b_1-m_2+1}=\sqbinom{n-k-1}{k}$ for some $0\leq k\leq n-1$,
i.e., $\Lambda\sim\binom{n-k+1,1,0}{k},\binom{k}{n-k+1,1,0}$ for some $0\leq k\leq n-1$.
\end{proof}

\subsection{Lusztig parametrizations for classical groups}\label{0234}
In this subsection, we recall the basic results of Lusztig's parametrization of
irreducible characters of finite classical groups.
First we define
\begin{align*}
\cals_{\rmO^+_{\rm even}} &=\{\,\Lambda\in\cals\mid{\rm def}(\Lambda)\equiv 0\pmod 4\,\}, &
\cals_{\rmO_{2n}^+} &=\{\,\Lambda\in\cals_{\rmO^+}\mid{\rm rk}(\Lambda)=n\,\}; \\
\cals_\Sp &=\{\,\Lambda\in\cals\mid{\rm def}(\Lambda)\equiv 1\pmod 4\,\}, &
\cals_{\Sp_{2n}} &=\{\,\Lambda\in\cals_\Sp\mid{\rm rk}(\Lambda)=n\,\}; \\
\cals_{\rmO^-_{\rm even}} &=\{\,\Lambda\in\cals\mid{\rm def}(\Lambda)\equiv 2\pmod 4\,\}, &
\cals_{\rmO_{2n}^-} &=\{\,\Lambda\in\cals_{\rmO^-}\mid{\rm rk}(\Lambda)=n\,\}.
\end{align*}
Then we have the following parametrization of irreducible unipotent characters from
\cite{lg} theorem 8.2:

\begin{prop}[Lusztig]\label{0232}
Let\/ $\bfG=\Sp_{2n}$ or $\rmO^\epsilon_{2n}$ for $\epsilon=\pm$.
Then there exists a bijective correspondence between $\cals_\bfG$ and $\cale(G)_1$.
\end{prop}

The following result (\cf.~\cite{DM} theorem 13.23, remark 13.24) is fundamental for the classification
of $\cale(G)$:

\begin{prop}[Lusztig]\label{0201}
Let\/ $\bfG=\Sp_{2n}$, $\SO_{2n+1}$, or $\rmO^\epsilon_{2n}$ for $\epsilon=\pm$,
and let $s\in(G^*)^0$.
There is a bijection $\grL_s\colon\cale(G)_s\rightarrow \cale(C_{G^*}(s))_1$
satisfying the condition
\[
\langle\rho,\epsilon_\bfG R^\bfG_{\bfT^*,s}\rangle_G
=\langle\grL_s(\rho),\epsilon_{C_{\bfG^*}(s)}R^{C_{\bfG^*}(s)}_{\bfT^*,{\bf1}}\rangle_{C_{G^*}(s)}
\]
for any rational maximal torus $\bfT^*$ containing $s$
where $C_{G^*}(s)$ denotes the centralizer of $s$ in\/ $G^*$,
$\epsilon_\bfG=(-1)^r$ and $r$ is the $\bfF_q$-rank of\/ $\bfG$,
and $\langle,\rangle_G$ denotes the inner product on the space of class functions on $G$.
\end{prop}
Such a bijection $\grL_s$ is called a \emph{Lusztig correspondence}.
The following proposition follows from \cite{lusztig-book} (9.9.1):

\begin{prop}[Lusztig]\label{0202}
The Lusztig series $\cale(G)_s$ has a cuspidal character
if and only if $\cale(C_{G^*}(s))_1$ has a (unipotent) cuspidal character
and the largest $\bfF_q$-split torus in the center of $C_{G^*}(s)$ coincides with
the largest $\bfF_q$-split torus in the center of\/ $G^*$.
In this case, $\rho$ is cuspidal if and only if $\grL_s(\rho)$ is cuspidal.
\end{prop}

\subsection{Theta correspondence for $(\Sp_{2n},\rmO^\epsilon_{2n'})$ and $(\Sp_{2n},\rmO_{2n'+1})$}\label{0208}

Let $\lambda=[\lambda_1,\ldots,\lambda_k]$ and $\mu=[\mu_1,\ldots,\mu_l]$ be two partitions with
$\lambda_1\geq\cdots\geq\lambda_k\geq 0$ and $\mu_1\geq\cdots\geq\mu_l\geq 0$.
We denote
\[
\lambda\preccurlyeq\mu\quad\text{ if \ }\mu_1\geq\lambda_1\geq\mu_2\geq\lambda_2\geq\cdots.
\]
Then we define two relations on the set $\cals$:
\begin{align}\label{0210}
\begin{split}
\calb^+ &=\{\,(\Lambda,\Lambda')\in\cals\times\cals
\mid\Upsilon(\Lambda_*)\preccurlyeq\Upsilon(\Lambda'^*),\ \Upsilon(\Lambda'_*)\preccurlyeq\Upsilon(\Lambda^*)\,\};\\
\calb^- &=\{\,(\Lambda,\Lambda')\in\cals\times\cals
\mid\Upsilon(\Lambda^*)\preccurlyeq\Upsilon(\Lambda'_*),\ \Upsilon(\Lambda'^*)\preccurlyeq\Upsilon(\Lambda_*)\,\}.
\end{split}
\end{align}
Then we define
\begin{align}\label{0227}
\begin{split}
\calb_{\Sp,\rmO^\epsilon_{\rm even}}
&=\{\,(\Lambda,\Lambda')\in\calb^+\cap(\cals_\Sp\times\cals_{\rmO^\epsilon_{\rm even}})\mid {\rm def}(\Lambda')
=-{\rm def}(\Lambda)+\epsilon\cdot 1\,\}, \\
\calb_{\Sp_{2n},\rmO_{2n'}^\epsilon}
&=\calb_{\Sp,\rmO^\epsilon_{\rm even}}\cap(\cals_{\Sp_{2n}}\times\cals_{\rmO_{2n'}^\epsilon})
\end{split}
\end{align}
where $\epsilon=\pm$.

Recall that $\omega_{\bfG,\bfG',1}$ denotes the unipotent part of Weil character
$\omega^\psi_{\bfG,\bfG'}$ for the dual pair $(\bfG,\bfG')=(\Sp_{2n},\rmO^\epsilon_{2n'})$.
The following proposition is from \cite{pan-finite-unipotent} theorem 3.4:
\begin{prop}\label{0220}
Let $(\bfG,\bfG')=(\Sp_{2n},\rmO^\epsilon_{2n'})$.
Then we have the decomposition
\[
\omega_{\bfG,\bfG',1}
=\sum_{(\Lambda,\Lambda')\in\calb_{\bfG,\bfG'}}\rho_\Lambda\otimes\rho_{\Lambda'}.
\]
\end{prop}

Let $\bfG$ be $\Sp_{2n}$, $\rmO^\epsilon_{2n}$ or $\SO_{2n+1}$.
For $s\in (G^*)^0$, we define
\begin{align}
\begin{split}
\bfG^{(0)}=\bfG^{(0)}(s)
&=\prod_{\langle\xi\rangle\subset\{\xi_1,\ldots,\xi_n\},\ \xi\neq\pm 1}\bfG_{[\xi]}(s); \\
\bfG^{(1)}=\bfG^{(1)}(s)
&=\bfG_{[-1]}(s); \\
\bfG^{(2)}=\bfG^{(2)}(s)
&=\begin{cases}\bfG_{[1]}(s), & \text{if $\bfG$ is an orthogonal group}; \\
(\bfG_{[1]}(s))^*, & \text{if $\bfG$ is a symplectic group}.
\end{cases}
\end{split}
\end{align}
where $\bfG_{[\xi]}(s)$ is given in \cite{amr} subsection 1.B
(see also \cite{pan-Lusztig-correspondence} subsection 2.2).
We know that $\bfG^{(0)}$ is a product of general linear groups or unitary group, and
\begin{equation}
(\bfG^{(1)},\bfG^{(2)})=\begin{cases}
(\Sp_{2n^{(1)}},\Sp_{2n^{(2)}}), & \text{if $\bfG=\SO_{2n+1}$};\\
(\rmO^{\epsilon^{(1)}}_{2n^{(1)}},\Sp_{2n^{(2)}}), & \text{if $\bfG=\Sp_{2n}$};\\
(\rmO^{\epsilon^{(1)}}_{2n^{(1)}},\rmO^{\epsilon^{(2)}}_{2n^{(2)}}), & \text{if $\bfG=\rmO^\epsilon_{2n}$}
\end{cases}
\end{equation}
for some non-negative integers $n^{(1)},n^{(2)}$ depending on $s$,
and some $\epsilon^{(1)},\epsilon^{(2)}$.
Note that if $\bfG=\rmO^\epsilon_{2n}$, then $\epsilon^{(1)},\epsilon^{(2)}$ also depend on $s$ (and $\epsilon$),
if $\bfG=\Sp_{2n}$, then $\epsilon^{(1)}$ can be $+$ or $-$ for each $s$ such that $n^{(1)}\geq 1$.
Now the element $s$ can be written as $s=s^{(0)}\times s^{(1)}\times s^{(2)}$
where $s^{(1)}$ (resp.~$s^{(2)}$) is the part whose eigenvalues are all equal to $-1$ (resp.~$1$),
and $s^{(0)}$ is the part whose eigenvalues do not contain $1$ or $-1$.
There is a natural bijection $\cale(C_{G^*}(s))_1\simeq\cale(G^{(0)}\times G^{(1)}\times G^{(2)})_1$,
and so from Proposition~\ref{0232}, we have a bijection
\begin{align}\label{0235}
\begin{split}
\Xi_s\colon\cale(G)_s &\rightarrow\cale(G^{(0)}\times G^{(1)}\times G^{(2)})_1 \\
\rho &\mapsto \rho^{(0)}\otimes\rho^{(1)}\otimes\rho^{(2)}
\end{split}
\end{align}
where $\rho^{(j)}\in\cale(G^{(j)})_1$, and $\Xi_s$ is called a \emph{modified Lusztig correspondence}.

Now suppose that $\bfG=\rmO_{2n+1}$.
Then $\bfG^0=\SO_{2n+1}$, $G\simeq G^0\times\{\pm1\}$, and we have
\[
\cale(G)=\{\,\rho\otimes{\bf 1},\rho\otimes\sgn\mid\rho\in\cale(G^0)\,\}.
\]
For $s\in(G^0)^*$ and $\epsilon=\pm$, we define
\begin{align*}
\cale(G)_{s,\epsilon}=\begin{cases}
\{\,\rho\otimes{\bf1}\mid\rho\in\cale(G^0)_{s}\,\}, & \text{if $\epsilon=+$};\\
\{\,\rho\otimes\sgn\mid\rho\in\cale(G^0)_{s}\,\}, & \text{if $\epsilon=-$}.
\end{cases}
\end{align*}
Each $\cale(G)_{s,\epsilon}$ is regarded as a Lusztig series for $\bfG=\rmO_{2n+1}$.
It is clear that
\[
\cale(G)=\bigcup_{(s)\subset(G^0)^*,\ \epsilon=\pm}\cale(G)_{s,\epsilon}.
\]

Suppose that $\bfG=\Sp_{2n}$ or $\rmO_{2n+1}$.
For $\rho\in\cale(G)$, we define $\epsilon_\rho=\pm1$ by $\rho(-1)=\epsilon_\rho\rho(1)$.
If $\bfG=\rmO_{2n+1}$ and $\rho\in\cale(G)_{s,\epsilon}$,
then it is obvious that $\epsilon_\rho=\epsilon$.

The following two propositions concerning the commutativity between the modified Lusztig correspondence
and the theta correspondence are from \cite{pan-Lusztig-correspondence} proposition~8.3 and proposition~8.8:

\begin{prop}\label{0209}
Let $(\bfG,\bfG')=(\Sp_{2n},\rmO_{2n'}^\epsilon)$ where $\epsilon=+$ or $-$,
and let $\rho\in\cale(G)_s$ and $\rho'\in\cale(G')_{s'}$ for some semisimple $s\in G^*$
and $s'\in (G'^*)^0$.
Write $\Xi_s(\rho)=\rho^{(0)}\otimes\rho^{(1)}\otimes\rho^{(2)}$ and\/
$\Xi_{s'}(\rho')=\rho'^{(0)}\otimes\rho'^{(1)}\otimes\rho'^{(2)}$
where $\Xi_s,\Xi_{s'}$ are the modified Lusztig correspondences for $\bfG$ and $\bfG'$ respectively
given in \cite{pan-ambiguity}.
Then $(\rho,\rho')$ occurs in $\Theta_{\bfG,\bfG'}^\psi$ if and only if
\begin{itemize}
\item $s^{(0)}=s'^{(0)}$ (up to conjugation), and $\rho^{(0)}=\rho'^{(0)}$;

\item $\bfG^{(1)}\simeq\bfG'^{(1)}$, and $\rho^{(1)}=\rho'^{(1)}$;

\item $(\rho^{(2)},\rho'^{(2)})$ occurs $\Theta_{\bfG^{(2)},\bfG'^{(2)},1}$
\end{itemize}
where $\Theta_{\bfG^{(2)},\bfG'^{(2)},1}$ denotes the unipotent part of\/ $\Theta^\psi_{\bfG^{(2)},\bfG'^{(2)}}$.
\end{prop}

\begin{prop}\label{0212}
Let $(\bfG,\bfG')=(\Sp_{2n},\rmO_{2n'+1})$,
and let $\rho\in\cale(G)_s$ and $\rho'\in\cale(G')_{s',\epsilon'}$ for some semisimple $s\in G^*$
and $s'\in G'^*$, and some $\epsilon'=\pm$.
Write $\Xi_s(\rho)=\rho^{(0)}\otimes\rho^{(1)}\otimes\rho^{(2)}$,
$\Xi_{s'}(\rho'|_{G'^0})=\rho'^{(0)}\otimes\rho'^{(1)}\otimes\rho'^{(2)}$
where $\Xi_s,\Xi_{s'}$ are the modified Lusztig correspondences for $\bfG$ and $\bfG'$
respectively given in \cite{pan-ambiguity}.
Then $(\rho,\rho')$ occurs in $\Theta_{\bfG,\bfG'}^\psi$ if and only if
\begin{itemize}
\item $s^{(0)}=-s'^{(0)}$ (up to conjugation), and $\rho^{(0)}=\rho'^{(0)}$;

\item $\bfG^{(2)}\simeq\bfG'^{(1)}$ and $\rho^{(2)}=\rho'^{(1)}$;

\item $(\rho^{(1)},\rho'^{(2)})$ occurs in $\Theta_{\bfG^{(1)},\bfG'^{(2)},1}$;

\item $\epsilon_\rho=\epsilon'$.
\end{itemize}
\end{prop}

\subsection{Theta correspondence for unitary dual pairs}\label{0233}
For a partition $\lambda$ of $n$ (we write $|\lambda|=n$),
let $\lambda_\infty$ denote its \emph{$2$-core}
and let $\Lambda_\lambda$ denote the symbol given in \cite{pan-Lusztig-correspondence} \S 5.
It is known that
\begin{equation}\label{0207}
|\lambda|=|\lambda_\infty|+2|\Upsilon(\Lambda_\lambda)|
\end{equation}
(\cf.~\cite{FS} p.233).
We define
\[
\cals_{\rmU_n}=\{\,\Lambda_\lambda\mid \lambda\text{ a partition of }n\,\}.
\]

Define $\calb^+_{\rmU,\rmU}$ to be the set consisting of pairs of symbols $(\Lambda,\Lambda')\in\calb^+$
(\cf.~(\ref{0210})) such that
\begin{equation}
{\rm def}(\Lambda')=\begin{cases}
-{\rm def}(\Lambda),-{\rm def}(\Lambda)+1, & \text{if ${\rm def}(\Lambda)$ is even};\\
-{\rm def}(\Lambda)+1,-{\rm def}(\Lambda)+2, & \text{if ${\rm def}(\Lambda)$ is odd},
\end{cases}
\end{equation}
Similarly, define $\calb^-_{\rmU,\rmU}$ to be the set consisting of pairs of symbols $(\Lambda,\Lambda')\in\calb^-$ such that
\begin{equation}
{\rm def}(\Lambda')=\begin{cases}
-{\rm def}(\Lambda)-2,-{\rm def}(\Lambda)-1, & \text{if ${\rm def}(\Lambda)$ is even};\\
-{\rm def}(\Lambda)-1,-{\rm def}(\Lambda), & \text{if ${\rm def}(\Lambda)$ is odd},
\end{cases}
\end{equation}
Then we define
\begin{equation}
\calb_{\rmU_n,\rmU_{n'}}=\{(\Lambda_\lambda,\Lambda_{\lambda'})\in\calb^+_{\rmU,\rmU}\cup\calb^-_{\rmU,\rmU}
\mid |\lambda|=n,\ |\lambda'|=n'\,\}.
\end{equation}
From \cite{pan-Lusztig-correspondence} lemma 5.1, it is known that
\begin{equation}\label{0214}
\calb_{\rmU_n,\rmU_{n'}}\subset\begin{cases}
\calb^+_{\rmU,\rmU}, & \text{if $n+n'$ is even};\\
\calb^-_{\rmU,\rmU}, & \text{if $n+n'$ is odd}.
\end{cases}
\end{equation}

The following proposition is rephrased from \cite{amr} th\'eor\`eme 5.15
(\cf.~\cite{pan-Lusztig-correspondence} proposition~5.13).

\begin{prop}\label{0519}
Let $(\bfG,\bfG')=(\rmU_n,\rmU_{n'})$.
Then we have the decomposition
\[
\omega_{\bfG,\bfG',1}
=\sum_{(\Lambda_\lambda,\Lambda_{\lambda'})\in\calb_{\bfG,\bfG'}}\rho_\lambda\otimes\rho_{\lambda'}
\]
where $\omega_{\bfG,\bfG',1}$ denotes the unipotent part of $\omega^\psi_{\bfG,\bfG'}$.
\end{prop}

For $s\in G^*$, we know that
$C_{\bfG^*}(s)=\prod_{\langle\xi\rangle}\bfG_{[\xi]}(s)$ where each $\bfG_{[\xi]}(s)$
is a unitary group or a general linear group.
Define
\begin{align*}
\bfG^{(0)} &=\prod_{\langle\xi\rangle,\ \xi\neq 1}\bfG_{[\xi]}(s)\\
\bfG^{(1)} &=\bfG_{[1]}(s)
\end{align*}
where $\bfG_{[\xi]}(s)$ is given as in Subsection~\ref{0208}.
Then $C_{\bfG^*}(s)=\bfG^{(0)}\times\bfG^{(1)}$ and hence there exists a (unique) one-to-one
correspondence
\begin{align}\label{0211}
\begin{split}
\Xi_s\colon \cale(G)_s &\rightarrow\cale(G^{(0)}\times G^{(1)})_1 \\
\rho &\mapsto \rho^{(0)}\otimes \rho^{(1)}
\end{split}
\end{align}
by Proposition~\ref{0232}.
Then $s$ can be written as $s=s^{(0)}\times s^{(1)}$
where $s^{(1)}$ is the part whose eigenvalues are all equal to $1$,
and $s^{(0)}$ is the part whose eigenvalues do not contain $1$.

The following proposition can be extracted from \cite{amr} th\'eor\`eme 2.6
(\cf.~\cite{pan-chain01} theorem~3.10).

\begin{prop}\label{0510}
Let $(\bfG,\bfG')=(\rmU_n,\rmU_{n'})$,
and let $\rho\in\cale(G)_s$ and $\rho'\in\cale(G')_{s'}$ for some semisimple $s,s'$.
Write $\Xi_s(\rho)=\rho^{(0)}\otimes\rho^{(1)}$ and\/ $\Xi_{s'}(\rho')=\rho'^{(0)}\otimes\rho'^{(1)}$.
Then $(\rho,\rho')$ occurs in $\Theta^\psi_{\bfG,\bfG'}$ if and only if
\begin{itemize}
\item $s^{(0)}=s'^{(0)}$, and $\rho^{(0)}=\rho'^{(0)}$,

\item $(\rho^{(1)},\rho'^{(1)})$ occurs in $\Theta_{\bfG^{(1)},\bfG'^{(1)},1}$.
\end{itemize}
\end{prop}

\section{Generalized Preservation Principle of Unipotent Characters}\label{0806}

\subsection{Preservation principle for unipotent cuspidal characters}\label{0206}
Irreducible unipotent cuspidal characters of classical groups (other than the general linear groups)
are characterized by Lusztig as follows:
\begin{enumerate}
\item $\rmU_n(q)$ has an irreducible unipotent character if
and only if $n=\frac{1}{2}m(m+1)$ for some $m\geq 0$.
In this case, there is a unique irreducible unipotent cuspidal character.

\item $\rmO_{2n+1}(q)$ has an irreducible unipotent character if
and only if $n=m(m+1)$ for some $m\geq 0$.
In this case, there are two irreducible unipotent cuspidal characters different by the sign character.

\item $\Sp_{2n}(q)$ has an irreducible unipotent character if
and only if $n=m(m+1)$ for some $m\geq 0$.
In this case, there is a unique irreducible unipotent cuspidal character.

\item $\rmO^\epsilon_{2n}(q)$ has an irreducible unipotent character if
and only if $n=m^2$ and $\epsilon=\sgn((-1)^m)$ for some $m\geq 0$.
In this case, there are two irreducible unipotent cuspidal characters different by the sign character.
\end{enumerate}

When $(\bfG,\bfG')=(\rmU_n,\rmU_{n'})$ or $(\Sp_{2n},\rmO^\epsilon_{2n'})$,
the unipotent characters are preserved by the Howe correspondence,
and the first occurrence of irreducible unipotent cuspidal characters are completely
determined in \cite{adams-moy} theorems 4.1 and 5.2.
And we have the following preservation principle of unipotent cuspidal characters:
\begin{enumerate}
\item[(I)] For the unipotent cuspidal character $\zeta^\rmU_m$ of $\rmU_{\frac{m(m+1)}{2}}(q)$,
we have
\[
n_0'^+(\zeta_m^\rmU)+n_0'^-(\zeta_m^\rmU)
=\tfrac{1}{2}m(m-1)+\tfrac{1}{2}(m+1)(m+2)=m(m+1)+1.
\]

\item[(II)] For the unipotent cuspidal characters $\zeta^{\rm I}_m,\zeta^{\rm II}_m$ of $\rmO^\epsilon_{2m^2}(q)$,
we have
\[
n_0'(\zeta_m^{\rm I})+n_0'(\zeta_m^{\rm II})=2m(m-1)+2m(m+1)=4m^2.
\]

\item[(III)] For the unipotent cuspidal character $\zeta_m^\Sp$ of $\Sp_{2m(m+1)}(q)$ and all the
orthogonal spaces $\bfv'^\pm$ are even-dimensional,
we have
\[
n_0'^+(\zeta_m^\Sp)+n_0'^-(\zeta_m^\Sp)=2m^2+2(m+1)^2=4m(m+1)+2.
\]
\end{enumerate}
Note that unipotent characters are not preserved by the Howe correspondence for the dual pair
$(\Sp_{2n},\rmO_{2n'+1})$.

\subsection{Generalized preservation principle for symplectic groups}
For a symbol $\Lambda=\binom{a_1,a_2,\ldots,a_{m_1}}{b_1,b_2,\ldots,b_{m_2}}\in\cals_\Sp$ and $\epsilon=\pm$,
we define $\theta^\epsilon_0(\Lambda)$ by
\begin{align*}
\theta_0^+(\Lambda) &=\begin{cases}
\binom{b_1,b_2,b_3,\ldots,b_{m_2}}{a_2,a_3,\ldots,a_{m_1}}=\Lambda^\rmt\smallsetminus\binom{-}{a_1},
& \text{if $m_1>0$};\\
\binom{b_1+1,b_2+1,\ldots,b_{m_2}+1,0}{-},
& \text{if $m_1=0$}.
\end{cases} \\
\theta_0^-(\Lambda) &=\begin{cases}
\binom{b_2,b_3,\ldots,b_{m_2}}{a_1,a_2,\ldots,a_{m_1}}=\Lambda^\rmt\smallsetminus\binom{b_1}{-},
& \text{if $m_2>0$};\\
\binom{-}{a_1+1,a_2+1,\ldots,a_{m_1}+1,0},
& \text{if $m_2=0$}.
\end{cases}
\end{align*}
It is known that ${\rm def}(\Lambda)\equiv 1\pmod 4$,
so it is easy to see that
${\rm def}(\theta_0^+(\Lambda))\equiv 0\pmod 4$ and ${\rm def}(\theta_0^-(\Lambda))\equiv 2\pmod 4$.
Then we have $\theta_0^\epsilon(\Lambda)\in\cals_{\rmO^\epsilon_{\rm even}}$.

\begin{lem}\label{0801}
For $\Lambda\in\cals_\Sp$ and $\epsilon=\pm$,
$\theta^\epsilon_0(\Lambda)$ is a symbol of minimal rank in
\[
\Theta^\epsilon(\Lambda)=\{\,\Lambda'\in\cals_{\rmO^\epsilon_{\rm even}}\mid(\Lambda,\Lambda')\in\calb^\epsilon\,\},
\]
i.e., $\theta^\epsilon_0(\Lambda)\in\Theta^\epsilon(\Lambda)$ and
${\rm rk}(\theta^\epsilon_0(\Lambda))\leq{\rm rk}(\Lambda')$ for any
$\Lambda'\in\Theta^\epsilon(\Lambda)$.
\end{lem}
\begin{proof}
From \cite{pan-uniform} lemma~2.6,
it is not difficult to see that $(\Lambda,\theta^\epsilon_0(\Lambda))\in\calb^\epsilon$
and $\theta^\epsilon_0(\Lambda)$ is of minimal rank in $\Theta^\epsilon(\Lambda)$.
\end{proof}

\begin{lem}\label{0802}
For $\Lambda\in\cals_\Sp$,
we have
\[
{\rm rk}(\theta_0^+(\Lambda))+{\rm rk}(\theta_0^-(\Lambda))
=2\,{\rm rk}(\Lambda)-\delta(\Lambda)+1.
\]
\end{lem}
\begin{proof}
Write $\Lambda=\binom{a_1,a_2,\ldots,a_{m_1}}{b_1,b_2,\ldots,b_{m_2}}\in\cals_\Sp$.
Note that $m_1-m_2\equiv 1\pmod 4$.
Now $m_1+m_2$ is odd and hence
\[
\bigl\lfloor\bigl(\tfrac{m_1+m_2}{2}\bigr)^2\bigr\rfloor-\bigl\lfloor\bigl(\tfrac{m_1+m_2-1}{2}\bigr)^2\bigr\rfloor
=\bigl\lfloor\bigl(\tfrac{m_1+m_2-1}{2}\bigr)^2\bigr\rfloor-\bigl\lfloor\bigl(\tfrac{m_1+m_2-2}{2}\bigr)^2\bigr\rfloor
= \tfrac{1}{2}(m_1+m_2-1).
\]
If $m_1=0$,
then $m_2\neq 0$ and
\begin{align*}
{\rm rk}(\theta_0^+(\Lambda)) &={\rm rk}(\Lambda)+m_2-\tfrac{1}{2}(m_1+m_2-1); \\
{\rm rk}(\theta_0^-(\Lambda)) &={\rm rk}(\Lambda)-b_1+\tfrac{1}{2}(m_1+m_2-1)
\end{align*}
by (\ref{0213}).
If $m_2=0$,
then $m_1\neq 0$ and
\begin{align*}
{\rm rk}(\theta_0^+(\Lambda)) &={\rm rk}(\Lambda)+m_1-\tfrac{1}{2}(m_1+m_2-1); \\
{\rm rk}(\theta_0^-(\Lambda)) &={\rm rk}(\Lambda)-a_1+\tfrac{1}{2}(m_1+m_2-1).
\end{align*}
If both $m_1,m_2$ are nonzero,
then
\begin{align*}
{\rm rk}(\theta_0^+(\Lambda)) &={\rm rk}(\Lambda)-a_1+\tfrac{1}{2}(m_1+m_2-1); \\
{\rm rk}(\theta_0^-(\Lambda)) &={\rm rk}(\Lambda)-b_1+\tfrac{1}{2}(m_1+m_2-1).
\end{align*}
For all three cases, we have
\[
{\rm rk}(\theta_0^+(\Lambda))+{\rm rk}(\theta_0^-(\Lambda))
=2\,{\rm rk}(\Lambda)-\delta(\Lambda)+1.
\]
\end{proof}

\begin{cor}\label{0805}
For an irreducible unipotent character $\rho_\Lambda$ of\/ $\Sp_{2n}(q)$,
we have
\[
n^+_0(\rho_\Lambda)+n_0^-(\rho_\Lambda)
=4n-2\delta(\Lambda)+2.
\]
\end{cor}
\begin{proof}
We know that ${\rm rk}(\Lambda)=n$ and
$n_0^\epsilon(\rho_\Lambda)=2\,{\rm rk}(\theta^\epsilon_0(\Lambda))$ by Lemma~\ref{0801}.
So the corollary follows from Lemma~\ref{0802} immediately.
\end{proof}

\subsection{Generalized preservation principle for even orthogonal groups}
For a symbol $\Lambda=\binom{a_1,a_2,\ldots,a_{m_1}}{b_1,b_2,\ldots,b_{m_2}}\in\cals_{\rmO^\epsilon_{\rm even}}$,
we define $\theta_0(\Lambda)$ by
\[
\theta_0(\Lambda)=\begin{cases}
\binom{b_1,b_2,b_3,\ldots,b_{m_2}}{a_2,a_3,\ldots,a_{m_1}}=\Lambda^\rmt\smallsetminus\binom{-}{a_1}, & \text{if $\epsilon=+$ and $m_1>0$};\\
\binom{b_1+1,b_2+1,\ldots,b_{m_2}+1,0}{-}, & \text{if $\epsilon=+$ and $m_1=0$}; \\
\binom{b_2,b_3,\ldots,b_{m_2}}{a_1,a_2,\ldots,a_{m_1}}=\Lambda^\rmt\smallsetminus\binom{b_1}{-}, & \text{if $\epsilon=-$ and $m_2>0$};\\
\binom{-}{a_1+1,a_2+1,\ldots,a_{m_1}+1,0}, & \text{if $\epsilon=-$ and $m_2=0$}.
\end{cases}
\]
If $\epsilon=+$, then ${\rm def}(\Lambda)\equiv 0\pmod 4$;
if $\epsilon=-$, then ${\rm def}(\Lambda)\equiv 2\pmod 4$.
Then we have ${\rm def}(\theta_0(\Lambda))\equiv 1\pmod 4$ for both cases,
and hence $\theta_0(\Lambda)\in\cals_\Sp$.

\begin{lem}\label{0804}
For $\Lambda\in\cals_{\rmO^\epsilon_{\rm even}}$,
$\theta_0(\Lambda)$ is a symbol of minimal rank in
\[
\Theta(\Lambda)=\{\,\Lambda'\in\cals_\Sp\mid(\Lambda',\Lambda)\in\calb^\epsilon\,\}.
\]
\end{lem}
\begin{proof}
The proof is similar to that of Lemma~\ref{0801}.
\end{proof}

\begin{lem}\label{0803}
For $\Lambda\in\cals_{\rmO^\epsilon_{\rm even}}$,
we have
\[
{\rm rk}(\theta_0(\Lambda))+{\rm rk}(\theta_0(\Lambda^\rmt))
=2\,{\rm rk}(\Lambda)-\delta(\Lambda).
\]
\end{lem}
\begin{proof}
Write $\Lambda=\binom{a_1,a_2,\ldots,a_{m_1}}{b_1,b_2,\ldots,b_{m_2}}$.
Note that now $m_1+m_2$ is even.
Suppose first that $m_1=m_2=0$.
Then $\Lambda=\Lambda^\rmt=\binom{-}{-}$ and $\epsilon=+$.
It is clearly that $\theta_0(\Lambda)=\theta_0(\Lambda^\rmt)=\binom{0}{-}$,
and the assertion is obvious.

Now we suppose that $m_1+m_2\geq 2$ and we have
\begin{align*}
\bigl\lfloor\bigl(\tfrac{m_1+m_2}{2}\bigr)^2\bigr\rfloor-\bigl\lfloor\bigl(\tfrac{m_1+m_2-1}{2}\bigr)^2\bigr\rfloor
&= \tfrac{1}{2}(m_1+m_2);\\
\bigl\lfloor\bigl(\tfrac{m_1+m_2-1}{2}\bigr)^2\bigr\rfloor-\bigl\lfloor\bigl(\tfrac{m_1+m_2-2}{2}\bigr)^2\bigr\rfloor
&= \tfrac{1}{2}(m_1+m_2-2).
\end{align*}
If $m_1=0$,
then $m_2\neq 0$ and
\[
{\rm rk}(\theta_0(\Lambda))+{\rm rk}(\theta_0(\Lambda^\rmt))
= {\rm rk}(\Lambda)+m_2-\tfrac{1}{2}(m_1+m_2)+{\rm rk}(\Lambda)-b_1+\tfrac{1}{2}(m_1+m_2-2).
\]
If $m_2=0$,
then $m_1\neq 0$ and
\[
{\rm rk}(\theta_0(\Lambda))+{\rm rk}(\theta_0(\Lambda^\rmt))
= {\rm rk}(\Lambda)+m_1-\tfrac{1}{2}(m_1+m_2)+{\rm rk}(\Lambda)-a_1+\tfrac{1}{2}(m_1+m_2-2).
\]
If both $m_1,m_2$ are nonzero,
then
\[
{\rm rk}(\theta_0(\Lambda))+{\rm rk}(\theta_0(\Lambda^\rmt))
=2\,{\rm rk}(\Lambda)-(a_1+b_1)+(m_1+m_2-2).
\]
For all cases we have
\[
{\rm rk}(\theta_0(\Lambda))+{\rm rk}(\theta_0(\Lambda^\rmt))=2\,{\rm rk}(\Lambda)-\delta(\Lambda).
\]
\end{proof}

\begin{cor}\label{0807}
For an irreducible unipotent character $\rho_{\Lambda}$ of\/ $\rmO^\epsilon_{2n}(q)$,
we have
\[
n'_0(\rho_\Lambda)+n'_0(\rho_{\Lambda}\cdot\sgn)
=4n-2\delta(\Lambda).
\]
\end{cor}
\begin{proof}
Now ${\rm rk}(\Lambda)={\rm rk}(\Lambda^\rmt)=n$,
$\rho_{\Lambda^\rmt}=\rho_\Lambda\cdot\sgn$, and $n'_0(\rho_\Lambda)=2\,{\rm rk}(\theta_0(\Lambda))$
by Lemma~\ref{0804}.
So the corollary follows from Lemma~\ref{0803} immediately.
\end{proof}

\subsection{Generalized preservation principle for unitary groups}
Suppose that $(\rho_\lambda,\rho_{\lambda'^\epsilon})$ occurs in the Howe correspondence for
the dual pair $(\rmU_n,\rmU_{n'^\epsilon})$ where $|\lambda|=n$, $|\lambda'^\epsilon|=n'^\epsilon$,
$\epsilon=\pm$ and $n'^+$ is an even integer, $n'^-$ is an odd integer.
\begin{enumerate}
\item Suppose that $n$ is even.
Then we know that $(\Lambda_\lambda,\Lambda_{\lambda'^+})\in\calb^+$ and $(\Lambda_\lambda,\Lambda_{\lambda'^-})\in\calb^-$
by (\ref{0214}).
To achieve the minimal rank, we see that
$\Lambda_{\lambda'^\epsilon}$ must be equal to $\theta^\epsilon_0(\Lambda_\lambda)$ given in Lemma~\ref{0801},
or equivalently, equal to $\theta_0(\Lambda_\lambda)$ with given $\epsilon$ in Lemma~\ref{0804}.

\item Suppose that $n$ is odd.
Then we know that $(\Lambda_\lambda,\Lambda_{\lambda'^-})\in\calb^+$ and $(\Lambda_\lambda,\Lambda_{\lambda'^+})\in\calb^-$.
To achieve the minimal rank, we see that
$\Lambda_{\lambda'^\epsilon}$ must be equal to $\theta^{-\epsilon}_0(\Lambda_\lambda)$ given in Lemma~\ref{0801},
or equivalently, equal to $\theta_0(\Lambda_\lambda)$ with given $-\epsilon$ in Lemma~\ref{0804}.
\end{enumerate}
Then from Proposition~\ref{0802} and Proposition~\ref{0803},
we know that
\begin{equation}\label{0810}
{\rm rk}(\theta_0^+(\Lambda_\lambda))+{\rm rk}(\theta_0^-(\Lambda_\lambda))
=\begin{cases}
2\,{\rm rk}(\Lambda_\lambda)-\delta(\Lambda_\lambda)+1, & \text{if ${\rm def}(\Lambda_\lambda)$ is odd};\\
2\,{\rm rk}(\Lambda_\lambda)-\delta(\Lambda_\lambda), & \text{if ${\rm def}(\Lambda_\lambda)$ is even}.
\end{cases}
\end{equation}

\begin{prop}\label{0811}
For an irreducible unipotent character $\rho_\lambda$ of\/ $\rmU_n(q)$,
we have
\[
n_0'^+(\rho_\lambda)+n_0'^-(\rho_\lambda)
=2n-2\delta(\rho_\lambda)+1.
\]
\end{prop}
\begin{proof}
Let $\lambda$ be a partition of $n$, and suppose $|\lambda_\infty|=\frac{1}{2}d(d+1)$ for some integer $d$.
If ${\rm def}(\Lambda_\lambda)$ is odd,
then by (\ref{0215}), (\ref{0207}) and (\ref{0810}), we conclude that
\begin{align*}
n_0^+(\rho_\lambda)+n_0^-(\rho_\lambda)
&= \frac{1}{2}d(d-1)+\frac{1}{2}(d+1)(d+2)+2\,{\rm rk}(\theta_0^+(\Lambda))+2\,{\rm rk}(\theta_0^-(\Lambda)) \\
&\qquad -2\biggl[\frac{-{\rm def}(\Lambda_\lambda)-1}{2}\biggr]^2-2\biggl[\frac{-{\rm def}(\Lambda_\lambda)+1}{2}\biggr]^2 \\
&= d(d+1)+1+4\,{\rm rk}(\Lambda_\lambda)
-4\biggl[\frac{{\rm def}(\Lambda_\lambda)-1}{2}\biggr]\biggl[\frac{{\rm def}(\Lambda_\lambda)+1}{2}\biggr]
-2\delta(\Lambda_\lambda) \\
&= 2\,|\lambda|-2\delta(\Lambda_\lambda)+1.
\end{align*}
If ${\rm def}(\Lambda_\lambda)$ is even,
then similarly we conclude that
\begin{align*}
 n_0^+(\rho_\lambda)+n_0^-(\rho_\lambda)
&= \frac{1}{2}d(d-1)+\frac{1}{2}(d+1)(d+2)+2\,{\rm rk}(\theta_0^+(\Lambda))+2\,{\rm rk}(\theta_0^-(\Lambda)) \\
&\quad -2\biggl[\frac{-{\rm def}(\Lambda_\lambda)-2}{2}\biggr]\biggl[\frac{-{\rm def}(\Lambda_\lambda)}{2}\biggr]
-2\biggl[\frac{-{\rm def}(\Lambda_\lambda)}{2}\biggr]\biggl[\frac{-{\rm def}(\Lambda_\lambda)+2}{2}\biggr] \\
&= d(d+1)+1+4\,{\rm rk}(\Lambda_\lambda) -4\biggl[\frac{{\rm def}(\Lambda_\lambda)}{2}\biggr]^2-2\delta(\Lambda_\lambda) \\
&= 2\,|\lambda|-2\delta(\Lambda_\lambda)+1.
\end{align*}
Note that $\delta(\rho_\lambda)$ is defined to be $\delta(\Lambda_\lambda)$.
\end{proof}

\section{General Generalized Preservation Principle}

Now we can deduce the generalized preservation principle for general irreducible
characters from Proposition~\ref{0510}, Proposition~\ref{0209}, Proposition~\ref{0212}
and the result for unipotent characters in Section~\ref{0806}.

\subsection{For unitary groups}
In this subsection, we consider type (I), i.e., $G=G(\bfv)$ is a unitary group.

\begin{thm}\label{0919}
Let $\rho$ be an irreducible character of\/ $\rmU(\bfv)$.
Then
\[
n_0'^+(\rho)+n_0'^-(\rho)=2\dim(\bfv)-2\delta(\rho)+1.
\]
\end{thm}
\begin{proof}
Let $G=\rmU(\bfv)$ and suppose that $\rho\in\cale(G)_s$ for some $s$ and write
$\Xi_s(\rho)=\rho^{(0)}\otimes\rho^{(1)}$ by (\ref{0211}).
Suppose that $\rho$ first occurs in the correspondences for both dual pairs
$(G(\bfv),G(\bfv'^+))$ and $(G(\bfv),G(\bfv'^-))$.
Recall that $C_{\bfG^*}=\prod_{\langle\xi\rangle}\bfG_{[\xi]}$ and
we identify $\bfG^*=\bfG$.
Then we have decompositions
\[
\bfv=\Biggl[\bigoplus_{\langle\xi\rangle,\ \xi\neq1}\bfv_{[\xi]}\Biggr]\oplus\bfv_{[1]},\qquad
\bfv'^\pm=\Biggl[\bigoplus_{\langle\xi\rangle,\ \xi\neq1}\bfv_{[\xi]}'^\pm\Biggr]\oplus\bfv_{[1]}'^\pm.
\]
By Proposition~\ref{0510} we have $\dim(\bfv_{[\xi]})=\dim(\bfv_{[\xi]}'^\pm)$ for $\xi\neq 1$.
Moreover the unipotent character $\rho^{(1)}$ first occurs in the correspondences
for both dual pairs $(G(\bfv_{[1]}),G(\bfv_{[1]}'^+))$ and $(G(\bfv_{[1]}),G(\bfv_{[1]}'^-))$.
Note that $\bfv_{[1]}'^+$ and $\bfv_{[1]}'^-$ are in different Witt series.
Hence we have
\[
\dim(\bfv'^+_{[1]})+\dim(\bfv'^-_{[1]})=2\dim(\bfv_{[1]})-2\delta(\rho^{(1)})+1
\]
by Proposition~\ref{0811}.
Therefore
\begin{align*}
\dim(\bfv'^+)+\dim(\bfv'^-)
& =\sum_{\langle\xi\rangle,\ \xi\neq1}\dim(\bfv_{[\xi]}'^+)+\dim(\bfv_{[1]}'^+)
+\sum_{\langle\xi\rangle,\ \xi\neq1}\dim(\bfv_{[\xi]}'^-)+\dim(\bfv_{[1]}'^-)\\
& =2\sum_{\langle\xi\rangle,\ \xi\neq1}\dim(\bfv_{[\xi]})
+ 2\dim(\bfv_{[1]})-2\delta(\rho^{(1)})+1 \\
& = 2\dim(\bfv)-2\delta(\rho)+1.
\end{align*}
Note that $\delta(\rho)=\delta(\rho^{(1)})$ by definition in (\ref{0107}).
\end{proof}

\subsection{For orthogonal groups}
In this subsection, we consider type (II), i.e., $G=G(\bfv)$ is an orthogonal group.

\begin{thm}\label{0920}
Let $\rho$ be an irreducible character of\/ $\rmO(\bfv)$.
Then
\[
n_0'(\rho)+n_0'(\rho\cdot\sgn)=2\dim(\bfv)-2\delta(\rho).
\]
\end{thm}
\begin{proof}
Let $G=\rmO(\bfv)$,
and suppose that $(\rho,\rho')$ (resp.~$(\rho\cdot\sgn,\rho'')$) first occurs in the correspondence
for the dual pair $(\rmO(\bfv),\Sp(\bfv'))$ (resp.~$(\rmO(\bfv),\Sp(\bfv''))$).

\begin{enumerate}
\item Suppose that $\bfv$ is even-dimensional and identify $G(\bfv)$ with its dual group.
Suppose that $\rho\in\cale(G)_s$ for some $s$ and write
$\Xi_s(\rho)=\rho^{(0)}\otimes\rho^{(1)}\otimes\rho^{(2)}$.
It is known that $\rho\cdot\sgn\in\cale(G)_s$ and
\[
\Xi_s(\rho\cdot\sgn)=\rho^{(0)}\otimes(\rho^{(1)}\cdot\sgn)\otimes(\rho^{(2)}\cdot\sgn).
\]
Now we have a decomposition
\[
\bfv=\Biggl[\bigoplus_{\langle\xi\rangle,\ \xi\neq\pm1}\bfv_{[\xi]}\Biggr]\oplus\bfv_{[-1]}\oplus\bfv_{[1]}
\]
such that
\[
G^{(0)}=\prod_{\langle\xi\rangle,\ \xi\neq\pm1}G(\bfv_{[\xi]}),\quad
G^{(1)}=G(\bfv_{[-1]}),\quad
G^{(2)}=G(\bfv_{[1]})^*=G(\bfv_{[1]}).
\]
Note that $G(\bfv_{[1]})$ is an even orthogonal group.

Let $\bfv'^*$ be the odd-dimensioanl quadratic space such that $G(\bfv'^*)$ is the dual group of $G(\bfv')$.
Then we have
\[
\bfv'^*
=\Biggl[\bigoplus_{\langle\xi\rangle,\ \xi\neq\pm1}\bfv_{[\xi]}'^*\Biggr]
\oplus\bfv_{[-1]}'^*\oplus\bfv_{[1]}'^*.
\]
Now $\bfv_{[1]}'^*$ is an odd-dimensional quadratic space.
Let $\bfv_{[1]}'$ be a symplectic space such that $G(\bfv_{[1]}')$ are the dual groups
of $G(\bfv_{[1]}'^*)$.
By Proposition~\ref{0209},
we have $\dim(\bfv_{[\xi]})=\dim(\bfv_{[\xi]}'^*)$ for $\xi\neq 1$.
Moreover, $\rho^{(2)}$ first occurs in the correspondences for the dual pairs
$(G(\bfv_{[1]}),G(\bfv_{[1]}'))$.
Similarly, we have a decomposition of $\bfv''^*$ as above,
and $\rho^{(2)}\cdot\sgn$ first occurs in the correspondences for the dual pairs
$(G(\bfv_{[1]}),G(\bfv_{[1]}''))$.
Hence by Corollary~\ref{0807},
we have
\[
\dim(\bfv_{[1]}')+\dim(\bfv_{[1]}'')=2\dim(\bfv_{[1]})-2\delta(\rho^{(2)}),
\]
and therefore
\begin{align*}
\dim(\bfv')+\dim(\bfv'')
& =\sum_{\langle\xi\rangle,\ \xi\neq1}\dim(\bfv'^*_{[\xi]})+\dim(\bfv'_{[1]})
+ \sum_{\langle\xi\rangle,\ \xi\neq1}\dim(\bfv''^*_{[\xi]})+\dim(\bfv''_{[1]}) \\
& =2\sum_{\langle\xi\rangle,\ \xi\neq1}\dim(\bfv_{[\xi]})+
2\dim(\bfv_{[1]})-2\delta(\rho^{(2)}) \\
& =2\dim(\bfv)-2\delta(\rho).
\end{align*}

\item Suppose that $\bfv$ is odd-dimensional, and $G=\rmO(\bfv)$, $G^0=\SO(\bfv)$.
Suppose that $\rho\in\cale(G)_{s,\epsilon}$ for some $s$ and some $\epsilon=\pm$.
Then $\rho\cdot\sgn\in\cale(G)_{s,-\epsilon}$,
and we can write
$\Xi_s(\rho|_{G^0})=\Xi_s(\rho\cdot\sgn|_{G^0})=\rho^{(0)}\otimes\rho^{(1)}\otimes\rho^{(2)}$.

Now $\bfv^*$ is an even-dimensional symplectic space with $\dim(\bfv^*)=\dim(\bfv)-1$
such that $G(\bfv^*)=G^0(\bfv)^*$,
and there is a decomposition
\[
\bfv^*=\Biggl[\bigoplus_{\langle\xi\rangle,\ \xi\neq\pm1}\bfv^*_{[\xi]}\Biggr]\oplus\bfv_{[1]}^* \oplus\bfv_{[-1]}^*
\]
where both $\bfv^*_{[1]}$ and $\bfv^*_{[-1]}$ are symplectic spaces,
and hence
\begin{equation}\label{0902}
\dim(\bfv)
=\dim(\bfv^*)+1
=\sum_{\langle\xi\rangle,\ \xi\neq -1}\dim(\bfv^*_{[\xi]})+\dim(\bfv_{[-1]}^*)+1.
\end{equation}

Suppose that $(\rho,\rho')$ first occurs in the correspondence for $(\rmO(\bfv),\Sp(\bfv'))$.
Then $\bfv'^*$ is odd-dimensional,
and there are decompositions
\[
\bfv'^*=\Biggl[\bigoplus_{\langle\xi\rangle,\ \xi\neq\pm1}\bfv'^*_{[\xi]}\Biggr]
\oplus\bfv_{[1]}'^* \oplus\bfv_{[-1]}'^*
\]
where $\bfv'^*_{[1]}$ is odd-dimensional and
$\bfv'^*_{[-1]}$ is even-dimensional.
Let $\bfv_{[1]}'$ be the symplectic space such that $G(\bfv_{[1]}')$ is the dual group of $G(\bfv_{[1]}'^*)$.
Then
\[
\dim(\bfv')
=\sum_{\langle\xi\rangle,\ \xi\neq 1}\dim(\bfv'^*_{[\xi]})+\dim(\bfv'_{[1]}).
\]
By Proposition~\ref{0212},
we have $\dim(\bfv^*_{[\xi]})=\dim(\bfv'^*_{[\xi]})$ for $\xi\neq \pm1$,
$\dim(\bfv_{[-1]}^*)=\dim(\bfv'_{[1]})$,
and $\rho^{(2)}$ first occurs in the correspondence for $(G(\bfv^*_{[1]}),G(\bfv'^*_{[-1]}))$.

Suppose that $(\rho\cdot\sgn,\rho'')$ first occurs in the correspondence for $(\rmO(\bfv),\Sp(\bfv''))$.
Similarly, we have $\dim(\bfv^*_{[\xi]})=\dim(\bfv''^*_{[\xi]})$ for $\xi\neq \pm1$,
$\dim(\bfv_{[-1]}^*)=\dim(\bfv''_{[1]})$,
and $\rho^{(2)}$ first occurs in the correspondence for $(G(\bfv^*_{[1]}),G(\bfv''^*_{[-1]}))$.

Now from \cite{pan-Lusztig-correspondence} lemma~8.6 and proposition~8.7, we known that
the two even-dimensional quadratic spaces $\bfv'^*_{[-1]}$ and $\bfv''^*_{[-1]}$ are in different Witt series.
Hence
\[
\dim(\bfv'^*_{[-1]})+\dim(\bfv''^*_{[-1]})
= 2\dim(\bfv^*_{[1]})-2\delta(\rho^{(2)})+2
\]
by Corollary~\ref{0805}.
Therefore by (\ref{0902}), we have
\begin{align*}
 \dim(\bfv')+\dim(\bfv'')
&= \sum_{\langle\xi\rangle,\ \xi\neq1}\dim(\bfv_{\xi}'^*)+\dim(\bfv'_{[1]})
+\sum_{\langle\xi\rangle,\ \xi\neq1}\dim(\bfv_{[\xi]}''^*)+\dim(\bfv''_{[1]}) \\
&=2\sum_{\langle\xi\rangle,\ \xi\neq\pm1}\dim(\bfv_{[\xi]}^*)+2\dim(\bfv^*_{[-1]})
+2\dim(\bfv^*_{[1]})-2\delta(\rho^{(2)})+2 \\
&=2\dim(\bfv)-2\delta(\rho).
\end{align*}
\end{enumerate}
Note that $\delta(\rho)=\delta(\rho^{(2)})$ for this situation.
\end{proof}

\begin{rem}
Suppose that $G(\bfv)=\rmO^+_{2n}(q)$, and let $\rho\in\cale(G)_s$ for some $s$.
Write $\Xi_s(\rho)=\rho^{(0)}\otimes\rho^{(1)}\otimes\rho^{(2)}$.
From (\ref{0107}) and Lemma~\ref{0216}, we know that
$\delta(\rho)=\delta(\rho^{(2)})\leq n$.
This means that if $n'+n''<n$, no irreducible character $\rho$ of $G$ such that
$\rho$ occurs in the correspondences for $(\rmO_{2n}^+,\Sp_{2n'})$ and
$\rho\cdot\sgn$ occurs in the correspondence for $(\rmO_{2n}^+,\Sp_{2n''})$.
Moreover, for $0\leq k\leq n$,
the unipotent character $\rho_{\binom{n-k}{k}}$ is the only irreducible character
$\rho$ of $G$ such that $\rho$ occurs in the correspondences for $(\rmO_{2n}^+,\Sp_{2k})$ and
$\rho\cdot\sgn$ occurs in the correspondence for $(\rmO_{2n}^+,\Sp_{2(n-k)})$.
\end{rem}

\begin{rem}
Suppose that $G(\bfv)=\rmO^-_{2n}(q)$, and let $\rho\in\cale(G)$ for some $s$.
Write $\Xi_s(\rho)=\rho^{(0)}\otimes\rho^{(1)}\otimes\rho^{(2)}$.
Now we know that $\delta(\rho)=\delta(\rho^{(2)})\leq n-1$.
This means that if $n'+n''<n+1$, no irreducible character $\rho$ of $G$ such that
$\rho$ occurs in the correspondences for $(\rmO_{2n}^-,\Sp_{2n'})$ and
$\rho\cdot\sgn$ occurs in the correspondence for $(\rmO_{2n}^-,\Sp_{2n''})$.
Moreover, by Lemma~\ref{0217} the unipotent character $\rho_{\binom{k}{n-k+1,1,0}}$ is the only irreducible character
$\rho$ of $G$ such that $\rho$ occurs in the correspondences for $(\rmO_{2n}^-,\Sp_{2k})$ and
$\rho\cdot\sgn$ occurs in the correspondence for $(\rmO_{2n}^-,\Sp_{2(n-k)+2})$.
\end{rem}

\subsection{For symplectic groups I}
In this subsection, we consider type (III), i.e., $G=G(\bfv)$ is a symplectic group,
and $G(\bfv'^\pm)$ are even orthogonal groups.

\begin{thm}\label{0921}
Let $\rho$ be an irreducible character of\/ $\Sp(\bfv)$.
Suppose that $\{\bfv'^\pm\}$ are Witt series of even-dimensional orthogonal spaces.
Then
\begin{equation}\label{0915}
n_0'^+(\rho)+n_0'^-(\rho)=2\dim(\bfv)-2\delta(\rho)+2.
\end{equation}
\end{thm}
\begin{proof}
Suppose that $(\rho,\rho'^+)$ (resp.~$(\rho,\rho'^-)$) first occurs in the correspondence
for the dual pair $(\Sp(\bfv),\rmO(\bfv'^+))$ (resp.~$(\Sp(\bfv),\rmO(\bfv'^-))$).
Suppose that $\rho$ is in $\cale(G)_s$ for some $s$ and
write $\Xi_s(\rho)=\rho^{(0)}\otimes\rho^{(1)}\otimes\rho^{(2)}$.
We know that $s$ determines a decomposition
\[
\bfv^*
=\Biggl[\bigoplus_{\langle\xi\rangle,\ \xi\neq\pm1}\bfv^*_{[\xi]}\Biggr]
\oplus\bfv^*_{[-1]}\oplus\bfv^*_{[1]}
\]
such that
\[
G^{(0)}=\prod_{\langle\xi\rangle,\ \xi\neq\pm1}G(\bfv^*_{[\xi]}),\quad
G^{(1)}=G(\bfv_{[-1]}^*),\quad
G^{(2)}=G(\bfv_{[1]}^*)^*.
\]
Let $\bfv_{[1]}$ be a $2\nu_1$-dimensional symplectic space such that $G(\bfv_{[1]})=G(\bfv_{[1]}^*)^*$.
Note that
\[
\dim(\bfv)
=\dim(\bfv^*)-1
=\sum_{\langle\xi\rangle,\ \xi\neq 1}\dim(\bfv_{[\xi]}^*)+\dim(\bfv_{[1]}).
\]
Let $\bfv_{[\xi]}'^{*\pm}$, $\bfv_{[-1]}'^{*\pm}$, and $\bfv_{[1]}'^{\pm}$ be defined similarly,
and hence
\[
\dim(\bfv'^\pm)=\sum_{\langle\xi\rangle,\ \xi\neq1}\dim(\bfv_{[\xi]}'^{*\pm})+\dim(\bfv_{[1]}'^\pm).
\]

By Proposition~\ref{0209}, we know that $\dim(\bfv^*_{[\xi]})=\dim(\bfv'^{*\pm}_{[\xi]})$ for $\xi\neq1$.
Moreover, $\rho^{(2)}$ first occurs in the correspondences for the dual pairs
$(G(\bfv_{[1]}),G(\bfv_{[1]}'^\pm))$.
Note that $\bfv_{[1]}'^+$ and $\bfv_{[1]}'^-$ are both even-dimensional quadratic spaces
and are in different Witt series.
Then by Corollary~\ref{0805},
we have
\[
\dim(\bfv'^+_{[1]})+\dim(\bfv'^-_{[1]})=2\dim(\bfv_{[1]})-2\delta(\rho^{(2)})+2.
\]
Therefore,
\begin{align*}
\dim(\bfv'^+)+\dim(\bfv'^-)
& =\sum_{\langle\xi\rangle,\ \xi\neq1}\dim(\bfv'^{*+}_{[\xi]})+\dim(\bfv'^+_{[1]})
+\sum_{\langle\xi\rangle,\ \xi\neq1}\dim(\bfv'^{*-}_{[\xi]})+\dim(\bfv'^-_{[1]}) \\
& =2\sum_{\langle\xi\rangle,\ \xi\neq1}\dim(\bfv_{[\xi]}^*)
+2\dim(\bfv_{[1]})-2\delta(\rho^{(2)})+2 \\
& = 2\dim(\bfv)-2\delta(\rho)+2.
\end{align*}
Note that now $\delta(\rho)=\delta(\rho^{(2)})$.
\end{proof}

\begin{rem}\label{0301}
Suppose that $G(\bfv)=\Sp_{2n}(q)$ and both $\bfv'^\pm$ are even-dimensional.
Let $\rho\in\cale(G)_s$ for some $s$,
and write $\Xi_s(\rho)=\rho^{(0)}\otimes\rho^{(1)}\otimes\rho^{(2)}$.
From (\ref{0107}) and Lemma~\ref{0216}, we know that
$\delta(\rho)=\delta(\rho^{(2)})\leq n$.
This means that if $n'+n''<n$, no irreducible character $\rho$ of $G$ such that
$\rho$ occurs in the correspondences for both $(\Sp_{2n},\rmO_{2n'}^+)$ and $(\Sp_{2n},\rmO^-_{2n''})$.
Moreover, if $\delta(\Lambda)=n$,
then $\Lambda\sim\Lambda_k:=\binom{n-k+1,0}{k}$
for some $0\leq k\leq n$.
Note that $\theta^+_0(\Lambda_k)=\binom{k}{0}\in\cals_{\rmO^+_{2k}}$ and
$\theta^-_0(\Lambda_k)=\binom{n-k+1,0}{-}\in\cals_{\rmO^-_{2(n-k)+2}}$.
Hence the unipotent character $\rho_{\Lambda_k}$ is the only irreducible character of $G$ which occurs in
the correspondences for both $(\Sp_{2n},\rmO^+_{2k})$ and $(\Sp_{2n},\rmO^-_{2(n-k)+2})$.
This is an refinement of \cite{howe-finite} corollary~5.2.1.
\end{rem}

\subsection{For symplectic groups II}\label{0918}
In this subsection, we consider type (III), i.e., $G=G(\bfv)$ is a symplectic group,
and $G(\bfv')$ is an odd orthogonal group.

\begin{thm}\label{0702}
Let $\rho$ be an irreducible character of\/ $\Sp(\bfv)$.
Suppose that $\{\bfv'\}$ is a Witt series of odd-dimensional orthogonal spaces.
Then
\begin{equation}\label{0916}
n_0'(\rho)+n_0'(\rho^c)=2\dim(\bfv)-2\delta(\rho)+2.
\end{equation}
\end{thm}
\begin{proof}
Let $\bfv$ be a symplectic space,
and let $\bfv^*$ be an odd-dimensional orthogonal space such that $G(\bfv^*)=G(\bfv)^*$.
We have a decomposition
\[
\bfv^*
=\Biggl[\bigoplus_{\langle\xi\rangle,\ \xi\neq\pm1}\bfv^*_{[\xi]}\Biggr]
\oplus\bfv^*_{[-1]}\oplus\bfv^*_{[1]}.
\]
Now $\bfv_{[1]}^*$ is a $(2\nu_1+1)$-dimensional orthogonal space.
Let $\bfv_{[1]}$ the $2\nu_1$-dimensional symplectic space such that $G(\bfv_{[1]})$ is the dual group
of $G(\bfv_{[1]}^*)$.
Then we have
\[
G^{(0)}=\prod_{\langle\xi\rangle,\ \xi\neq\pm1}G(\bfv^*_{[\xi]}),\qquad
G^{(1)}=G(\bfv_{[-1]}^*),\qquad
G^{(2)}=G(\bfv_{[1]}).
\]
So now
\[
\dim(\bfv)=\dim(\bfv^*)-1
=\sum_{\langle\xi\rangle,\ \xi\neq 1}\dim(\bfv^*_{[\xi]})+\dim(\bfv_{[1]}),
\]
and $\bfv^*_{[-1]}$ is an even-dimensional orthogonal space.

Also $\bfv'$ is an odd-dimensional orthogonal space,
then $\bfv'^{*}$ is a symplectic space,
and we have a decomposition
\[
\bfv'^{*}
=\Biggl[\bigoplus_{\langle\xi\rangle,\ \xi\neq\pm1}\bfv'^*_{[\xi]}\Biggr]
\oplus\bfv'^*_{[-1]}\oplus\bfv'^*_{[1]}.
\]
Therefore, we have
\[
\dim(\bfv')
=\sum_{\langle\xi\rangle,\ \xi\neq 1}\dim(\bfv'^{*}_{[\xi]})+\dim(\bfv'^{*}_{[1]})+1,
\]
and $\bfv'^{*}_{[-1]}$ is symplectic space.

Let $\rho$ be an irreducible character in $\cale(G)_s$ for some $s$ and write
$\Xi_s(\rho)=\rho^{(0)}\otimes\rho^{(1)}\otimes\rho^{(2)}$.
From \cite{pan-ambiguity} \S 8, we know $\rho^c\in\cale(G)_s$ and
\begin{equation}\label{0410}
\Xi_s(\rho^c)=\rho^{(0)}\otimes(\rho^{(1)}\cdot\sgn)\otimes\rho^{(2)}.
\end{equation}
Now suppose that $\rho$ (resp.~$\rho^c$) first occurs in the correspondence for the pair
$(G(\bfv),G(\bfv'))$ (resp.~$(G(\bfv),G(\bfv''))$) where $\bfv',\bfv''$ are in the same Witt series
of odd-dimensional quadratic spaces.
By Proposition~\ref{0212},
we have that $\dim(\bfv^*_{[\xi]})=\dim(\bfv'^{*}_{[\xi]})=\dim(\bfv''^{*}_{[\xi]})$ for
$\xi\neq 1$.
Moreover, $\rho^{(1)}$ (resp.~$\rho^{(1)}\cdot\sgn$) first occurs in the correspondence for the dual pair
$(G(\bfv^*_{[-1]}),G(\bfv'^{*}_{[-1]}))$ (resp.~$(G(\bfv^*_{[-1]}),G(\bfv''^{*}_{[-1]}))$).
Hence
\[
\dim(\bfv'^{*}_{[-1]})+\dim(\bfv''^{*}_{[-1]})
= 2\dim(\bfv^*_{[-1]})-2\delta(\rho^{(1)})
\]
by Corollary~\ref{0807}.
Therefore
\begin{align*}
\dim(\bfv')+\dim(\bfv'')
&= \sum_{\langle\xi\rangle,\ \xi\neq 1}\dim(\bfv_{[\xi]}'^{*})+\dim(\bfv_{[1]}'^{*})+1
+\sum_{\langle\xi\rangle,\ \xi\neq1}\dim(\bfv_{[\xi]}''^{*})+\dim(\bfv_{[1]}''^{*})+1\\
&=2\sum_{\langle\xi\rangle,\ \xi\neq\pm1}\dim(\bfv_{[\xi]}^*)+2\dim(\bfv_{[1]})+2
+2\dim(\bfv^*_{[-1]})-2\delta(\rho^{(1)}) \\
&=2\dim(\bfv)-2\delta(\rho)+2.
\end{align*}
Note that now $\delta(\rho^c)=\delta(\rho^{(1)}\cdot\sgn)=\delta(\rho^{(1)})=\delta(\rho)$.
\end{proof}

Let $\{\bfv'^\pm\}$ be the two Witt series of odd-dimensional quadratic spaces.
For $\rho\in\cale(\Sp(\bfv))$,
let $n_0'^\pm(\rho)$ denote the minimal dimension of $\bfv'^\pm$ such that
$\rho$ occurs in the correspondence for the dual pair $(G(\bfv),G(\bfv'^\pm))$.

\begin{lem}\label{0903}
Let $\rho\in\cale(\Sp(\bfv))$.
Then $n_0'^-(\rho^c)=n_0'^+(\rho)$.
\end{lem}
\begin{proof}
Let $\psi$ be a non-trivial additive character of $\bfF_q$.
It is well known that $\omega^{\psi_a}_{\Sp(\bfv),\rmO(\bfv'^+)}=\omega^\psi_{\Sp(\bfv),\rmO(\bfv'^-)}$
where $\psi_a$ is another additive character of $\bfF_q$ given by
$\psi_a(x)=\psi(ax)$ and $a$ is a non-square element in $\bfF_q^\times$.
Then by \cite{pan-ambiguity} lemma~5.2, we know that
$\rho$ occurs in the correspondence for the dual pair $(\Sp(\bfv),\rmO(\bfv'^+))$ if and only if
$\rho^c$ occurs in the correspondence for the dual pair $(\Sp(\bfv),\rmO(\bfv'^-))$.
Hence the lemma is proved.
\end{proof}

Now we have another formulation of Theorem~\ref{0702}:
\begin{cor}
Let $\rho$ be an irreducible character of\/ $\Sp(\bfv)$.
Suppose that $\{\bfv'^\pm\}$ are the two Witt series of odd-dimensional quadratic spaces.
Then
\[
n_0'^+(\rho)+n_0'^-(\rho)=2\dim(\bfv)-2\delta(\rho)+2.
\]
\end{cor}

\begin{rem}
Suppose that $G(\bfv)=\Sp_{2n}(q)$ and $\bfv'$ is odd-dimensional.
Let $\rho\in\cale(G)_s$ for some $s$,
and write $\Xi_s(\rho)=\rho^{(0)}\otimes\rho^{(1)}\otimes\rho^{(2)}$.
From (\ref{0107}) and Lemma~\ref{0216}, we know that
$\delta(\rho)=\delta(\rho^{(1)})\leq n$.
This means that if $n'+n''<n$, no irreducible character $\rho$ of $G$ such that
$\rho$ occurs in the correspondences for $(\Sp_{2n},\rmO_{2n'+1})$ and
$\rho^c$ occurs in the correspondences for $(\Sp_{2n},\rmO_{2n''+1})$.

Recall that an irreducible character $\rho$ of $G$ is pseudo-unipotent
if $\rho\in\cale(G)_s$ such that both $G^{(0)}(s),G^{(2)}(s)$ are trivial.
In this case, we have $G^{(1)}(s)=\rmO^\epsilon_{2n}(q)$ where $\epsilon=\pm$ (\cf.~Lemma~\ref{0408}).
Then by (\ref{0107}) and Lemma~\ref{0215}, we see that $\delta(\rho)=n$ if and only if
$\rho\in\cale(G)_s$ such that $G^{(1)}(s)=\rmO^+_{2n}(q)$ and $\rho^{(1)}=\rho_{\binom{n-k}{k}}$
for some $0\leq k\leq n$.
Denote this pseudo-unipotent character by $\rho_k$.
Now $\rho_k$ is the only irreducible character $\rho$ of $G$ such that
$\rho$ occurs in the correspondences for $(\Sp_{2n},\rmO_{2k+1})$ and
$\rho^c$ occurs in the correspondences for $(\Sp_{2n},\rmO_{2(n-k)+1})$.
\end{rem}

\section{Pseudo-unipotent Cuspidal Characters}

\subsection{Pseudo-unipotent cuspidal characters}
Pseudo-unipotent characters are defined in Subsection~\ref{0101}.
By convention,
we regard the trivial character of the trivial group $\Sp_0(q)=\rmO_0^+(q)=\SO_1(q)$ to be both
unipotent and pseudo-unipotent.

Recall that for a general linear group, a unitary group, a special orthogonal group or an orthogonal group
$\bfG$, a character $\chi_\bfG$ of order $2$ is given in \cite{pan-odd} \S 3.3.

\begin{lem}\label{0405}
Let $\bfG$ be a unitary group, a special orthogonal group or an orthogonal group.
Then $\rho\in\cale(G)$ is pseudo-unipotent if and only if $\rho\chi_\bfG$ is unipotent, i.e.,
$\rho$ is pseudo-unipotent if and only if it belongs to the Lusztig series $\cale(G)_{-1}$.
\end{lem}
\begin{proof}
For such a classical group $\bfG$,
$-1$ is in the center of $G^*$,
and hence is in every rational maximal torus $\bfT^*$ of $\bfG^*$.
It is known that $R_{\bfT_w,\theta_w}=\chi_\bfG R_{\bfT_w,1}$ for any $w\in\bfW_\bfG$
(\cf.~\cite{DM} proposition 13.30),
and the mapping
\[
\cale(G)_{-1}\rightarrow\cale(G)_{1}\quad\text{given by }\rho\mapsto\rho\chi_\bfG
\]
is a bijection.
Therefore, $\langle\rho, R_{\bfT_w,\theta_w}\rangle_G\neq0$ if and only if
$\langle\rho\chi_\bfG, R_{\bfT_w,1}\rangle_G\neq0$.
\end{proof}

From the previous lemma we see that a pseudo-unipotent cuspidal character of a unitary group or 
an orthogonal group $G$ is exactly of the form $\zeta\chi_G$ where 
$\zeta$ is a unipotent cuspidal character of $G$.

\begin{lem}\label{0408}
Let $\bfG=\Sp_{2n}$.
Then $\rho\in\cale(G)$ is pseudo-unipotent if and only if it is
in the Lusztig series $\cale(G)_s$ where $s\in G^*$ such that
$C_{\bfG^*}(s)\simeq\rmO_{2n}^\epsilon$ where $\epsilon=\pm$.
\end{lem}
\begin{proof}
Let $\rho$ be an irreducible pseudo-unipotent character of $G=\Sp_{2n}(q)$.
Then we have $\langle\rho,R_{\bfT_w,\theta_w}\rangle_G\neq 0$ for some $(\bfT_w,\theta_w)$.
Let $(\bfT^*,s)$ be the pair corresponding to $(\bfT_w,\theta_w)$ where
$\bfT^*$ is a rational maximal torus of $\SO_{2n+1}$ and $s\in T^*$.
From the definition of $\theta_w$ we know that each $\lambda_i$ of
$s=(\lambda_1,\ldots,\lambda_n)\in T^*$ is $-1$.
Moreover, $\bfT^*$ is either a maximal torus of $\rmO_{2n}^+\subset\SO_{2n+1}$
or a maximal torus of $\rmO_{2n}^-\subset\SO_{2n+1}$.
Therefore, $C_{\bfG^*}(s)\simeq\rmO_{2n}^\epsilon$ for $\epsilon=\pm$.
\end{proof}

\begin{proof}[Proof of Proposition~\ref{0404}]
Let $\bfG=\Sp_{2n}$.
Consider the semisimple $s=s^\epsilon$ in $G^*$ such that $C_{\bfG^*}(s)\simeq\rmO_{2n}^\epsilon$.
By Proposition~\ref{0202}
for a pseudo-unipotent character $\rho\in\cale(\Sp_{2n}(q))_s$ to be cuspidal,
the unipotent character $\grL_s(\rho)$ of $\rmO^\epsilon_{2n}(q)$ must be cuspidal,
i.e., we must have $n=m^2$ and $\epsilon=(-1)^m$ for some integer $m$.
When $m>0$, from Subsection~\ref{0206},
we know that $\rmO_{2m^2}^{(-1)^m}(q)$ has two unipotent cuspidal characters.
Therefore $\Sp_{2m^2}(q)$ has exactly two pseudo-unipotent cuspidal characters.
\end{proof}

\begin{rem}
Let $\rho_1,\rho_2$ denote the two pseudo-unipotent cuspidal characters of $\Sp_{2m^2}(q)$
given by Proposition~\ref{0404}.
Because the two unipotent cuspidal characters of $\rmO^{(-1)^m}_{2m^2}$ are different by the sign character,
we have $\rho_2=\rho_1^c$ (\cf.~(\ref{0410})).
\end{rem}

Next we want to investigate all the first occurrences $(\rho,\rho')$ of unipotent or pseudo-unipotent
cuspidal characters for symplectic-orthogonal dual pairs.

\subsection{Correspondence for $(\Sp_{2n},\rmO^\epsilon_{2n'})$}
Suppose that $\rho$ is a unipotent cuspidal character of\/ $\Sp_{2m(m+1)}(q)$,
and $(\rho,\rho')$ is a first occurrence of cuspidal characters for
the dual pair $(\Sp_{2m(m+1)},\rmO_{2n'}^\epsilon)$.
Then from \cite{adams-moy} theorem 5.2 we know that either
\begin{enumerate}
\item $n'=m^2$ and $\epsilon=(-1)^m$; or

\item $n'=(m+1)^2$ and $\epsilon=(-1)^{m+1}$.
\end{enumerate}
Moreover, $\rho'$ is a unipotent cuspidal character of $\rmO_{2n'}^\epsilon(q)$.

Now we consider the first occurrence of pseudo-unipotent cuspidal characters.

\begin{prop}\label{0409}
Suppose that $\rho$ is a pseudo-unipotent cuspidal character of\/ $\Sp_{2m^2}(q)$
and $(\rho,\rho')$ is a first occurrence of cuspidal characters for $(\Sp_{2m^2},\rmO_{2n'}^\epsilon)$.
Then $n'=m^2$ or $m^2+1$ (depending on $\epsilon$).
Moreover, if $n'=m^2$,
$\rho'$ is also pseudo-unipotent.
\end{prop}
\begin{proof}
Let $(\bfG,\bfG')=(\Sp_{2m^2},\rmO^\epsilon_{2n'})$.
Suppose $\rho\in\cale(G)_s$ for some $s$ and write $\Xi_s(\rho)=\rho^{(0)}\otimes\rho^{(1)}\otimes\rho^{(2)}$
and $\rho'\in\cale(G')_{s'}$ and $\Xi_{s'}(\rho')=\rho'^{(0)}\otimes\rho'^{(1)}\otimes\rho'^{(2)}$.
By Lemma~\ref{0408} and the definition in Subsection~\ref{0208},
we know that $\bfG^{(0)}$ is trivial, $\bfG^{(1)}=\rmO_{2m^2}^{\epsilon'}$ and $\epsilon'=(-1)^m$,
$\bfG^{(2)}=\Sp_0$,
and $\rho^{(1)}$ is a unipotent cuspidal character of $\rmO_{2m^2}^{\epsilon'}(q)$.

By Proposition~\ref{0209} we know that both $\bfG'^{(0)}$ is trivial,
$\bfG'^{(1)}=\rmO_{2m^2}^{\epsilon'}$,
$\rho^{(1)}=\rho'^{(1)}$,
and $(\rho^{(2)},\rho'^{(2)})$ is a first occurrence of unipotent cuspidal characters for the dual pair
$(\bfG^{(2)},\bfG'^{(2)})$.
Then by \cite{adams-moy} theorem 5.3, we know that $\bfG'^{(2)}$ is either $\rmO_{0}^+$ or $\rmO_{2}^-$.
\begin{enumerate}
\item Suppose $\bfG'^{(2)}=\rmO_{0}^+$.
Then $s'$ is the element $-1\in G'^*=\rmO_{2m^2}^\epsilon(q)$ with $\epsilon=\epsilon'=(-1)^m$.
Hence $n'=m^2$ and $\rho'$ is pseudo-unipotent.

\item Suppose $\bfG'^{(2)}=\rmO_{2}^-$.
Then $s'$ is the element $(-1,1)\in C_{\bfG^*}(s')=\rmO_{2m^2}^{\epsilon'}\times\rmO^-_2\subset \bfG^*$
with $\epsilon=-\epsilon'=(-1)^{m+1}$.
Hence $n'=m^2+1$ and $\rho'$ is not pseudo-unipotent (and nor unipotent).
\end{enumerate}
\end{proof}

\subsection{Correspondence for $(\Sp_{2n},\rmO_{2n'+1})$}

\begin{prop}\label{0704}
Suppose $\rho$ is the unipotent cuspidal character of\/ $\Sp_{2m(m+1)}(q)$,
and $(\rho,\rho')$ is a first occurrence for $(\Sp_{2m(m+1)},\rmO_{2n'+1})$.
Then $n'=m(m+1)$ and $\rho'$ is a pseudo-unipotent cuspidal character of\/ $\rmO_{2m(m+1)+1}(q)$.
\end{prop}
\begin{proof}
Let $(\bfG,\bfG')=(\Sp_{2m^2},\rmO_{2n'+1})$, $\bfG'^0=\SO_{2n'+1}$, $\rho'^0=\rho'|_{G'^0}$.
Then $\rho'^0$ is an irreducible character.
As in the proof of Proposition~\ref{0409}, suppose that $\rho$ is in $\cale(G)_s$
and write $\Xi_s(\rho)=\rho^{(0)}\otimes\rho^{(1)}\otimes\rho^{(2)}$,
$\rho'_0\in\cale(G'^0)_{s'}$, $\Xi_{s'}(\rho'^0)=\rho'^{(0)}\otimes\rho'^{(1)}\otimes\rho'^{(2)}$,
$C_{\bfG^*}(s)=\rmO_{2m^2}^{\epsilon}$, and $\epsilon=(-1)^m$.
Now both $\bfG^{(0)}$ and $\bfG^{(1)}$ are trivial,
$\bfG^{(2)}=\Sp_{2m(m+1)}$,
and $\rho^{(2)}=\rho$.
By Proposition~\ref{0212} we know that $\bfG'^{(0)}$ is trivial,
$\bfG'^{(1)}\simeq \bfG^{(2)}=\Sp_{2m(m+1)}$,
and $(\rho^{(1)},\rho'^{(2)})$ is a first occurrence
of unipotent cuspidal characters for the dual pair $(\bfG^{(1)},\bfG'^{(2)})$.
This means that $\bfG'^{(2)}=\Sp_0$ and $\rho'^{(2)}$ is the trivial character.
Therefore, $n'=m(m+1)$ and $\rho'$ is a pseudo-unipotent cuspidal character of\/ $\rmO_{2m(m+1)+1}(q)$.
\end{proof}

\begin{prop}\label{0407}
Suppose that $\rho$ is a pseudo-unipotent cuspidal character of\/ $\Sp_{2m^2}(q)$,
and $(\rho,\rho')$ is a first occurrence of cuspidal characters for $(\Sp_{2m^2},\rmO_{2n'+1})$.
Then $n'=m(m-1)$ or $m(m+1)$,
and $\rho'$ is a unipotent cuspidal character of\/ $\rmO_{2n'+1}(q)$.
\end{prop}
\begin{proof}
Let $(\bfG,\bfG')=(\Sp_{2m^2},\rmO_{2n'+1})$, $\bfG'^0=\SO_{2n'+1}$, $\rho'^0=\rho'|_{G'^0}$.
As in the proof of Proposition~\ref{0704}, suppose that $\rho$ is in $\cale(G)_s$
$\rho'^0\in\cale(G'^0)_{s'}$, 
and write $\Xi_s(\rho)=\rho^{(0)}\otimes\rho^{(1)}\otimes\rho^{(2)}$,
$\Xi_{s'}(\rho'^0)=\rho'^{(0)}\otimes\rho'^{(1)}\otimes\rho'^{(2)}$,
$C_{\bfG^*}(s)=\rmO_{2m^2}^{\epsilon}$, and $\epsilon=(-1)^m$.
Now both $\bfG^{(0)}$ and $\bfG^{(2)}$ are trivial,
$\bfG^{(1)}=\rmO_{2m^2}^\epsilon$,
and $\rho^{(1)}$ is a unipotent cuspidal character of $\rmO_{2m^2}^\epsilon(q)$.
By Proposition~\ref{0212} we know that both $\bfG'^{(0)}$ and $\bfG'^{(1)}$ are trivial,
and $(\rho^{(1)},\rho'^{(2)})$ is a first occurrence
of unipotent cuspidal characters for the dual pair $(\bfG^{(1)},\bfG'^{(2)})$.
Then by \cite{adams-moy} theorem 5.3, we know that $\bfG'^{(2)}$ is $\Sp_{2m(m-1)}$ or $\Sp_{2m(m+1)}$.
We see that $s'$ is the identity element in $(\bfG'^0)^*=\Sp_{2n'}$
and this implies that $n'$ is $m(m-1)$ or $m(m+1)$, and $\rho'$ is unipotent.
\end{proof}

\bibliography{refer}
\bibliographystyle{amsalpha}

\end{document}